\documentclass[12pt]{article}
\usepackage[russian]{babel}
\usepackage{graphicx}
\usepackage[bf]{caption2}
\usepackage{amsfonts, amssymb, amsthm, amsmath}
\usepackage{color}
\newcounter{theorem}

\newcounter{theoremcounter}

\newcounter{remarkcounter}

\newcounter{corollarycounter}

\newtheorem{theorem}[theoremcounter]{Theorem}

\newtheorem{remark}[remarkcounter]{Remark}

\newtheorem{corollary}[corollarycounter]{Corollary}
\newcommand{\la}{\lambda}
\newcommand{\vp}{{\mathbf p}}
\newcommand{\vz}{{\mathbf z}}
\newcommand{\vf}{{\mathbf f}}
\newcommand{\vx}{{\mathbf x}}

\begin{document}

\title{Two approaches to the construction of perturbation bounds for continuous-time
Markov chains}

\author{Alexander Zeifman\thanks{Vologda State
University, Vologda, Russia; Institute of Informatics Problems,
Federal Research Center ``Computer Science and Control'',  Russian
Academy of Sciences, Moscow, Russia; Vologda Research Center of the
Russian Academy of Sciences, Vologda, Russia. E-mail:
a$\_$zeifman@mail.ru} \and  Victor Korolev\thanks{Faculty of
Computational Mathematics and Cybernetics, Moscow State University,
Moscow, Russia; Hanghzhou Dianzi University, Hangzhou, China;
Institute of Informatics Problems, Federal Research Center
``Computer Science and Control'', Russian Academy of Sciences,
Moscow, Russia} \and Yacov Satin\thanks{Vologda State University,
Vologda, Russia}}

\date{}

\maketitle

\sloppy

{\bf Abstract.} The paper is largely of a review nature. It
considers two main methods used to study stability and obtain
appropriate quantitative estimates of perturbations of
(inhomogeneous) Markov chains with continuous time and a finite or
countable state space. An approach is described to the construction
of perturbation estimates for the main five classes of such chains
associated with queuing models. Several specific models are
considered for which the limit characteristics and perturbation
bounds for admissible "perturbed" processes are calculated.

\bigskip


\section{Introduction}

In the paper, some topics are considered that  are related with the
stability of both homogeneous and non-homogeneous continuous-time
Markov chains with respect to the perturbation of their intensities
(infinitesimal characteristics). It is assumed that the evolution of
the system under consideration is described by a Markov chain with
the known state space and it is the  infinitesimal
matrix that is given inexactly. Different classes of admissible
perturbations can be considered. The ``perturbed'' infinitesimal
matrix can be arbitrary, and the small deviation of its norm from
that of the original matrix is assumed or it can be assumed that the
structure of the  infinitesimal matrix is known and
only its elements are ``perturbed'' within the same structure. Below
we will give a detailed description of these cases. In some papers
it is assumed that the perturbations have a special form, for
example, are expanded in a power series of a small parameter. This
assumption seems to be too restrictive and unrealistic.

The study of stability of characteristics of  stochastic models has
been actively developing since the 1970-s \cite{kalash,sto,zol}. At
that time V. M. Zolotarev proposed an approach within which limit
theorems of probability theory were treated as stability theorems.
This approach was cultivated in the works of V. M. Zolotarev, V. V.
Kalashnikov, V. M. Kruglov, V. Yu. Korolev and their colleagues in
the framework of international seminars on stability problems for
stochastic models founded by V. M. Zolotarev (see the series of the
proceedings of the seminar published as Springer Lecture Notes
starting from \cite{kz1} or as issues of the Journal of Mathematical
Sciences). This approach proved to be very productive for the study
of random sums in queueing theory, renewal theory and theory of
branching processes \cite{gk}.

Since 1980-s the problems related to the  estimation of stability of
Markov chains with respect to perturbations of their characteristics
have been thoroughly studied by N. V. Kartashov for homogeneous
discrete-time chains with general state space and, in parallel, by
A. I. Zeifman for inhomogeneous continuous-time chains within the
seminar mentioned above, see \cite{kar81,kar85,z85}. In particular,
a general approach for inhomogeneous continuous-time chains was
developed in \cite{z85}. In that paper both uniform and strong cases
were considered. Later birth-death processes were considered in
\cite{z88} and general properties and estimates for inhomogeneous
finite chains were considered in \cite{z94}. The paper \cite{z98}
was specially devoted to estimates for general birth-death processes
with the queueing system $M_t|M_t|N$ considered as an example. It
should be mentioned that these papers were not noticed in western
papers. For example, in \cite{aan} it was stated that there were NO
papers on stability of the (simplest stationary!!) system $M|M|1$.
For the first time we used the term `perturbation bounds' instead of
`stability' in the paper \cite{z12} on the referee's prompt. The
same situation takes place with the Kartashov's papers cited above.
The methods proposed in those papers seem to be used by most authors
of subsequent studies in estimation of perturbations of
discrete-time chains. Possibly, poor acquaintance with the early
papers of Kartashov and Zeifman can be explained by the differences
in terminology mentioned above: in the original (and foundational)
papers the term 'stability' was used (in the proceedings of the
seminar with the consonant appellation '{\it Stability} Problems for
Stochastic Models').

The present paper deals only with continuous-time chains, therefore the subsequent remarks mainly regard such a case.

Note that to obtain explicit and exact estimates of perturbation bounds of a chain, it is required to have estimates of the rate of convergence of the chain to its limit characteristics in the form of explicit inequalities. Moreover, the sharper convergence rate estimates, the more accurate perturbation bounds. These bounds can be more easily obtained for finite homogeneous Markov chains. Therefore, most publications concern just this situation, see, e. g., \cite{af,ds,mit03,mit04,mit05,mit06}. As this is so, two main approaches can be highlighted.

The first of them can be used for the case of weak ergodicity of a
chain in the uniform operator topology. The first bounds in this
direction were obtained in \cite{z85}. The principal progress
related to the replacement of the constant $S$ with $\log S$ in the
bound was implemented in \cite{mit03} and continued in Mitrophanov's
papers \cite{mit04,mit05,mit06} for the case of homogeneous chains
and then in \cite{Zeifman2011,z12} and in the subsequent papers of
these authors for the inhomogeneous chains. The contemporary state
of affairs in this field and new applied problems related to the
link between convergence rate and perturbation bounds in the
`uniform' case were described in \cite{mitrophanov2018}.
 In some recent papers uniform perturbation bounds of
homogeneous Markov chains were studied by the techniques of
stochastic differential equations, see for instance \cite{Shao} and
the references therein.

The second approach is used in the case where the uniform
ergodicity is not assured, that is typical for the processes most
interesting from the practical viewpoint. For example, birth-death
processes used for modeling queueing systems, and real processes in
biology, chemistry, physics, as a rule, are {\it not} uniformly
ergodic.

Following the ideas of N. V. Kartashov  (see a detailed description
in \cite{kar96}), most authors use the probability methods to study
ergodicity and perturbation bounds of stationary chains (with a
finite, countable or general state space) in various norms
\cite{aan,fhl,ma}. For a wide class of (mainly) stationary
discrete-time chains a close approach was considered in \cite{mt}
and more recent papers
\cite{abbas,jiang,liu2012,liu2018,medina,negrea,rudolf,thiede,truquet,vial,zheng}.

In the works of the authors of the present paper perturbation bounds for non-stationary finite or infinite continuous-time chains were studied by other methods.

The first papers dealing with non-statonary queueing models appeared
in the 1970-s (see \cite{g,g1}, and the more recent paper
\cite{mw}). Moreover, as far back as in \cite{gm} it was noted that
it is principally possible to use the logarithmic matrix norm for
the study of convergence rate of continuous-time Markov chains. The
corresponding general approach employing the theory of differential
equations in Banach spaces was developed in a series of papers by
the authors of the present paper, see a detailed description in
\cite{Zeifman2014inf,Zeifman2014s}. In \cite{z85} (also see
\cite{z88,z94}) a method for the study of perturbation bounds for
the vector of state probabilities of a continuous-time Markov chain
with respect to the perturbations of infinitesimal characteristics
of the chain in the total variation norm  ($l_1$-norm) was proposed.
The paper \cite{z98} contained a detailed study of the stability of
essentially non-stationary birth-death processes with respect to
conditionally small perturbations. Convergence rate estimates in
terms of weight norms and hence, the corresponding bounds for new
classes of Markov chains were considered in
\cite{Zeifman2014q,Zeifman2018c,Zeifman2019mat2,Zeifman2019amc}.

In the present paper both approaches are considered as well as the classes of inhomogeneous Markov chains for which at least one of these approaches yields reasonable perturbation bounds for basic probability characteristics.

The paper is organized as follows. In Section 2 basic notions and preliminary results are introduced. In Section 3 general theorems on perturbation bounds are considered. Section 4 contains convergence rate estimates and perturbation bounds for basic classes of the chains under consideration. Finally, in Section 5 some special queueing models are studied.

\section{Basic notions and preliminaries}

Let $X=X(t)$, $t\geq 0$, be, in general, inhomogeneous continuous-time Markov chain with a finite or countable state space $E_S= 0,1,\dots,S$, $S \le \infty$. The transition probabilities for $X=X(t)$ will be denoted $p_{ij}(s,t)=\Pr\left\{ X(t)=j\left| X(s)=i\right. \right\}$, $i,j
\ge 0, \;0\leq s\leq t$. Let $p_i(t)=\Pr\left\{ X(t) =i \right\}$ be the state probabilities of the chain and ${\bf p}(t) = \left(p_0(t), p_1(t), \dots\right)^T$ be the corresponding vector of state probabilities. In what follows it is assumed that
$$
\Pr\left\{X\left( t+h\right) =j|X\left( t\right) =i\right\} = \nonumber
$$
\begin{equation}\label{100}
=
\begin{cases}
q_{ij}\left( t\right)  h+\alpha_{ij}\left(t, h\right), & \mbox { if }j\neq i\vspace{1mm}
\cr 1-{\displaystyle\sum\limits_{k\neq i}q_{ik}\left( t\right)  h+\alpha_{i}\left(
t,h\right)}, & \mbox { if } j=i,
\end{cases}
\end{equation}
where all $\alpha_{i}(t,h)$ are $o(h)$ uniformly in $i$, that is, $\sup_i |\alpha_i(t,h)| = o(h)$.

As usual, we assume that if a chain is inhomogeneous, then all the
 infinitesimal characteristics (intensity functions)
 $q_{ij}\left( t\right)$ are  integrable in $t$ on any
 interval $[a,b]$, $0 \le a \le b$.

Let $a_{ij}(t) =  q_{ji}(t)$ for $j\neq i$ and $a_{ii}(t) =
-\sum_{j\neq i} a_{ji}(t) = -\sum_{j\neq i} q_{ij}(t)$.

Further, to provide the possibility to obtain more evident estimates
we will assume that 
\begin{equation}
|a_{ii}(t)| \le L < \infty
\label{0102-1}
\end{equation}
\noindent for almost all $t \ge 0$.

Then the state probabilities satisfy the forward Kolmogorov system
\begin{equation} \label{ur01}
\frac{d\vp}{dt}=A(t)\vp(t),
\end{equation}
\noindent where  $A(t)=Q^T(t)$, and $Q(t)$ is the
infinitesimal matrix of the process.

 Let  $\|\cdot\|$ be the usual  $l_1$-norm, i.e.
$\|{\vx}\|=\sum|x_i|$, and $\|B\| = \sup_j \sum_i |b_{ij}|$ for $B =
(b_{ij})_{i,j=0}^{\infty}$.  Denote $\Omega =
\left\{{\vx}\,:\, {\vx} \in l_1^+\  \& \ \|{\vx}\|=1\right\}$.
 Then
$$\|A(t)\| = 2\sup_{k}\left|a_{kk}(t)\right| \le 2 L $$
for almost all $t \ge 0$, and we can apply all results of \cite{DK} to equation (\ref{ur01})
in the space $l_1$.

Let $p_0(t) = 1 - \sum_{i \ge 1} p_i(t)$. Then from
(\ref{ur01}) we obtain the following equation (for a detailed
discussion, see, e. g., \cite{gz04,z06}):
\begin{equation}
\frac{d\vz}{dt}= B(t)\vz(t)+\vf(t), \label{2.06}
\end{equation}
\noindent where $\vf(t)=\left(a_{10}, a_{20},\cdots \right)^T$,
\begin{equation}
{
B = \left(
\begin{array}{ccccccccc}
a_{11}- a_{10}   & a_{12} - a_{10}   &  \cdots & a_{1r} - a_{10}  &\cdots \\
a_{21} - a_{20} & a_{22} - a_{20}   &    \cdots & a_{2r} - a_{20} &\cdots \\
a_{31} - a_{30}    & a_{32} - a_{30}  &    \cdots & a_{3r} - a_{30}  &\cdots \\
\cdots \\
a_{r1} - a_{r0}  & a_{r2} - a_{r0} & \cdots     &  a_{rr} -a_{r0} &\cdots \\
\cdots \cdots \cdots
\end{array}
\right).}\label{2.07}
\end{equation}

By $\bar{X}=\bar{X}(t)$ we will denote the `perturbed' Markov chain
with the same state space, state probabilities $\bar{p}_i(t)$,
transposed infinitesimal matrix $\bar{A}(t) =
\left(\bar{a}_{ij}(t)\right)_{i,j=0}^{\infty}$ and so on, and the
`perturbations' themselves, that is, the differences between the
corresponding `perturbed' and original characteristics will be
denoted by $\hat{a}_{ij}(t)$, $\hat{A}(t)$.

 Let $E(t,k) =
E\left\{X(t)\left|X(0)=k\right.\right\}$. Recall that a Markov chain
$X(t)$ is weakly ergodic, if $\|{\bf p^{*}}(t) - {\bf p^{**}}(t)\|
\to 0$ as $t \to \infty$ for any initial condition, and it has the
 limiting mean $\phi(t)$, if  $|E(t,k) -
\phi(t) | \to 0$ as $t \to \infty$ for any $k$.

Now briefly describe the main classes of the chains under consideration. The details concerning the first four classes can be found in
\cite{Zeifman2018lncs,Zeifman2018c}.

{\bf Class I.} Let $a_{ij}(t)=0$ for all $t\ge 0$ è $|i-j|>1$,
$a_{i,i+1}(t)=\mu_{i+1}(t)$,  $a_{i+1,i}(t)=\la_i(t)$. This is an inhomogeneous birth-death process (BDP) with the intensities $\la_i(t)$ (of birth) and $\mu_{i+1}(t)$ (of death) correspondingly.

{\bf Class II.} Now let $a_{ij}(t)=0$ for $i<j-1$,
$a_{i+k,i}(t)=a_k(t)$ for $k\ge 1$, and $a_{i,i+1}(t)=\mu_{i+1}(t)$.
This chain describes, for instance, the number of customers in a
queueing system in which the customers arrive in groups, but are
served one by one (in this case $a_k(t)$ is the arrival intensity of
a group of $k$ customers, and $\mu_{i}(t)$ is the service intensity
of the $i$th customer). The simplest models of this type were
considered in \cite{Nelson1988}, also see
\cite{Zeifman2018lncs,Zeifman2018c}.

{\bf Class III.} Let $a_{ij}(t)=0$ for $i>j+1$,
$a_{i,i+k}(t)=b_k(t)$, $k\ge 1$, and $a_{i+1,i}(t)=\la_i(t)$.  This
situation occurs in modelling queueing systems with arrivals of
single customers and group service.

{\bf Class IV.} Let $a_{i+k,i}(t)=a_k(t)$, $a_{i,i+k}(t)=b_k(t)$ for ${k\ge 1}$. This process appears in the description of a system with group arrival and group service, for earlier studies see \cite{s11,s12,Zeifman2014q}.

{\bf Class V.} Consider a Markov chain with `catastrophes' used for
modelling of some queueing systems, see, e. g., \cite{du1,du2,Di
Crescenzo2008,Li2017,Zhang2015,z12}.  Here the intensities have a
general form whereas a single (although substantial) restriction
consists in that the zero state is attainable from any other state
and the corresponding intensities $q_{k,0}(t)= a_{0,k}(t)$ for $
k\ge 1$ are called the intensities of catastrophes.

Now consider the following example illustrating some specific features of the problem under consideration.

{\bf Example} \cite{z12}. Consider a homogeneous BDP (class I) with
the intensities $\lambda_k(t) = 1$, $\mu_k(t) = 4$ for all $t$, $k$
and denote by $A$ the corresponding transposed intensity matrix.
Then, as is known (see, e. g., \cite{z95}), the BDP is strongly
ergodic and stable in the corresponding norm. On the other hand,
take a perturbed process with the transposed infinitesimal matrix
$\bar{A}=A+\hat{A}$, where $\hat{a}_{00}=-\varepsilon$,
$\hat{a}_{k0}=\frac{\varepsilon}{k(k+1)}$ for $k \ge 1$ and
$\hat{a}_{ij}= 0$ for the other $i,j$. Then the perturbed Markov
chain $\bar{X}(t)$ (describing the 'M|M|c queue with mass arrivals
when empty', see \cite{Chen1997,Li2017,Zhang2015}) is not ergodic,
since from the condition $\bar{A}{\bf \bar{p}} =0$ it follows that
the coordinates of the stationary distribution (if it exists) must
satisfy the condition $4\bar{p}_{k+1}= \bar{p}_{k} +
\bar{p}_{0}\frac{\varepsilon}{k+1} \ge
\bar{p}_{0}\frac{\varepsilon}{k+1}$, which is impossible.

As it has already been noted,  the (upper) bounds of perturbations
are closely connected with the (correspondingly, upper) estimates
for the convergence rate (also see the two next sections). On the
other side, it is also possible to construct important lower
estimates of the rate of convergence providing that the influence of
the initial conditions cannot fade too rapidly, see
\cite{Zeifman2018spl}. It turns out that it is principally
impossible to construct lower bounds for perturbations. Indeed, if
we consider the same BDP and as a perturbed BDP choose a BDP with
the intensities $\bar{\lambda}_k(t) = 1+\varepsilon$,
$\bar{\mu}_k(t) = 4(1+\varepsilon)$, then the stationary
distribution for the perturbed process will be the same as for the
original BDP for any positive $\varepsilon$.

\section{General theorems concerning perturbation bounds}

First consider {\it uniform} bounds.

\begin{theorem}
Let the Markov chain $X(t)$ be exponentially weakly ergodic, that is, for any initial conditions ${\bf p^{*}}(s) \in \Omega$,
${\bf p^{**}}(s) \in \Omega$ and any $s \ge 0$, $t \ge s$ there holds the inequality
\begin{equation}
\|{\bf p^{*}}(t) - {\bf p^{**}}(t)\|  \le 2c e^{-b(t-s)}.
\label{207}
\end{equation}
Then, for the perturbations small enough ($\hat{A}(t) \le \varepsilon$ for almost all $t \ge 0$), the perturbed chain $\bar{X}(t)$ is also exponentially weakly ergodic and the following perturbation bound takes place:
\begin{equation}
\limsup_{t \to \infty}   \|{\bf p}(t)- {\bar{\bf p}}(t)\| \le
\frac{\left(1+ \log (c/2)\right)\varepsilon}{b} . \label{207''}
\end{equation}
\end{theorem}

For the {\it proof} we will use the approach proposed in \cite{mit03} and modified in \cite{Zeifman2011} for the inhomogeneous case, also see \cite{z12}. Let \begin{equation}
\beta (t, s)=\sup_{ \| {\bf v} \| =1, \sum {v_i}=0} {\|U(t){\bf
v}\|} =
{\textstyle \frac{1}{2}} \sup_{i,j} \sum_k {|p_{ik}(t, s)-p_{jk}(t,
s)|}.
\label{208}
\end{equation}
Then
\begin{equation}
\|{\bf p}(t)-{\bf \bar{p}}(t)\| \le \beta(t,s)\|{\bf p}(s)-{\bf
\bar{p}}(s)\| + \int_s^t \|\hat{A}(u)\| \beta(u, s) du, \label{209}
\end{equation}
\noindent moreover,
\begin{equation}
\beta (t,s) \le 1, \quad \beta (t,s) \le {\textstyle\frac{c}{2}} e^{-b(t-s)}, \ 0
\le s \le t. \label{210}
\end{equation}

$$
\|{\bf p}(t)-\bar{\bf p}(t)\| \le
$$
\begin{equation}
=\begin{cases}
\|{\bf p}(s)-{\bf \bar{p}}(s)\|+ (t-s)\varepsilon, &  0<t\!-\!s<{\textstyle\frac1b} \ln {\frac{c}{2}},\vspace{1mm} \cr
{\frac{c}{2}e^{-b(t-s)} \|{\bf p}(s)-{\bf \bar{p}}(s)\|\!+\!\textstyle\frac1b}(\ln {\frac{c}{2}}\! +1\!-\!ce^{-b(t-s)})\varepsilon, &  t\!-\!s \ge
{\textstyle\frac1b}\ln {\frac{c}{2}},\end{cases} \label{211}
\end{equation}

\noindent whence, as $t \to \infty$, we obtain (\ref{207''}).

\begin{corollary} If under the conditions of theorem 1 the Markov chain $X(t)$ has a finite state space, then both Markov chains $X(t)$ and $\bar{X}(t)$ have limit expectations and
\begin{equation}
|\phi (t) - \bar{\phi} (t)| \le{\textstyle\frac1b}  S\left(1+ \log
(c/2)\right)\varepsilon. \label{207'''} \end{equation}
\end{corollary}

Now turn to the {\it weighted} bounds. Here we use the approach proposed in \cite{z85}, also see the detailed description in \cite{Zeifman2014inf,Zeifman2014s}.

 Let $1 \le d_1 \le
d_2 \le \dots$,
\begin{equation}
D=\left(
\begin{array}{ccccccc}
d_1   & d_1 & d_1 & \cdots  \\
0   & d_2  & d_2  &   \cdots  \\
0   & 0  & d_3  &   \cdots  \\
& \ddots & \ddots & \ddots \\
\end{array}
\right). \label{2013} \end{equation} \noindent Let  $l_{1D}=\left\{{\bf z} = (p_1,p_2,\cdots)^T
\,:\, \|{\bf z}\|_{1D} \equiv \|D {\bf z}\| <\infty \right\}$. Then
$\|B\|_{1D} = \|DBD^{-1}\|$. In addition, let $\|{\bf p}\|_{1D} =
\|{\bf z}\|_{1D}$.

Below we will assume that the following conditions hold:
\begin{equation}
\|B(t)\|_{1D} \le  \mathfrak{B}< \infty, \quad  \|{\bf f}(t)\|_{1D} \le \mathfrak{f}< \infty
\label{3001} \end{equation}
\noindent for almost all $t \ge 0$.

\medskip

Recall that  $X(t)$ is  $1D$-exponentially weakly ergodic Markov
chain, if
\begin{equation}
\|{\bf p^{*}}(t) - {\bf p^{**}}(t)\|_{1D} \le M
e^{-a\left(t-s\right)} \|{\bf p^{*}}(s) - {\bf p^{**}}(s)\|_{1D}.
\label{d-erg2}
\end{equation}
\noindent  for some   $M>0$,  $a>0$ and  any $s, t$:
$t\ge s \ge 0$,  any initial conditions ${\bf p^{*}}(s) \in l_{1D}$,
${\bf p^{**}}(s) \in l_{1D}$.

\smallskip
 If one can choose ${\bf p^{**}}(t) ={\bf \pi}$, then
the chain is $1D$-exponentially strongly ergodic.

Let
\begin{equation}
\|B(t)-\bar{B}(t)\|_{1D} \le  \left|\mathfrak{B-\bar{B}}\right|,
\quad \|{\bf f}(t)-\bar{{\bf f}}(t)\|_{1D} \le
\left|\mathfrak{f-\bar{f}}\right|. \label{3001'}
\end{equation}
for almost all $t \ge 0$.

\begin{theorem} \label{nonh-pert1}
If a Markov chain $X(t)$ is $1D$-exponentially weakly ergodic, then
$\bar{X}(t)$ is also $1D$-exponentially weakly ergodic and the
following perturbation estimate in the $1D$-norm holds:
\begin{eqnarray}
\limsup_{t \to \infty}   \|{\bf p}(t)- {\bar{\bf p}}(t)\|_{1D} \le
\frac{M\left(M\left|\mathfrak{B-\bar{B}}\right| \mathfrak{f}  +
a\left|\mathfrak{f-\bar{f}}\right|\right)}{a\left(a -
M\left|\mathfrak{B-\bar{B}}\right|\right)}.
\label{3002}
\end{eqnarray}
\noindent If $W=\inf_{i \ge 1} \frac {d_i}{i} > 0$, then both chains
$X(t)$ and $\bar{X}(t)$ have limiting means and
\begin{eqnarray}
\limsup_{t \to \infty}  |\phi(t) - \bar{\phi}(t)|\le
\frac{M\left(M\left|\mathfrak{B-\bar{B}}\right| \mathfrak{f}  +
a\left|\mathfrak{f-\bar{f}}\right|\right)}{Wa\left(a -
M\left|\mathfrak{B-\bar{B}}\right|\right)}.
\label{3003}\end{eqnarray}
\end{theorem}

\renewcommand{\proofname}{Proof.}

\begin{proof} The detailed consideration can be found in \cite{Zeifman2014s}.
Here we only outline the scheme of reasoning. Let  $V(t, s)$ and
$\bar{V}(t, s)$ be  the Cauchy operators for equation (\ref{2.06}) and for
the corresponding 'perturbed' equation, respectively. Then
\begin{equation}
\|V(t,s)\|_{1D} \le   M e^{-a(t-s)}, \quad  \|\bar{V}(t,s)\|_{1D}
\le   M e^{-\left(a -
M\left|\mathfrak{B-\bar{B}}\right|\right)(t-s)}  \label{2015}
\end{equation}
for all $t \ge s \ge 0$.  Then, rewriting  equation (\ref{2.06}) as
\begin{equation}
\frac{d\vz}{dt}=\bar{B}(t)\vz(t) + {\bf
f}(t)+\left(B(t)-\bar{B}(t)\right)\vz(t), \label{2020}
\end{equation}
\noindent after some algebra we obtain the following inequality in
the $1D$-norm:
$$
\left\|\vz(t)-\bar{\vz}(t)\right\|\le \int_0^t \|\bar{V}(t,\tau)\|
\left(\| B (\tau)-\bar{B} (\tau)\| \|\vz(\tau)\| + \|
\vf(\tau) -\bar{\vf}(\tau)\|\right)\ d\tau \le
$$
\begin{equation}
\le \int_0^t M e^{-\left(a -
M\left|\mathfrak{B-\bar{B}}\right|\right)(t-\tau)}
\left(\left|\mathfrak{B-\bar{B}}\right| \|\vz(\tau)\| +
\left|\mathfrak{f-\bar{f}}\right|\right)\ d\tau. \label{3005}
\end{equation}
On the other hand, $ \|\vz (t)\|_{1D} \le  Me^{-at}\|\vz (0)\|_{1D}
+ \frac{M}{a}  \mathfrak{f}$, for any $0 \le s \le t$.
Hence, under any initial condition $\vz (0) \in l_{1D}$ we obtain the following inequalities for the $1D$-norm:
$$
\left\|\vz(t)-\bar{\vz}(t)\right\| \le M\left(\left|\mathfrak{B-\bar{B}}\right|\frac{M}{a}
 \mathfrak{f}  + \left|\mathfrak{f-\bar{f}}\right|\right) \int_0^t e^{-\left(a -
 M\left|\mathfrak{B-\bar{B}}\right|\right)(t-\tau)}\, d\tau +
$$
$$
+ M \int_0^t e^{-\left(a -
M\left|\mathfrak{B-\bar{B}}\right|\right)(t-\tau)}
\left|\mathfrak{B-\bar{B}}\right| Me^{-a\tau}\|\vz (0)\| \, d\tau \le
$$
\begin{equation}
\le \frac{M\left(M\left|\mathfrak{B-\bar{B}}\right| \mathfrak{f} +
a\left|\mathfrak{f-\bar{f}}\right|\right)}{a\left(a -
M\left|\mathfrak{B-\bar{B}}\right|\right)} + o\left( 1 \right).\label{2021-1}
\end{equation}
\noindent So, the first assertion of the theorem is proved.

 Then second assertion follows from the inequality
$\|\vz\|_{1E}\le W^{-1}\|\vz\|_{1D}$ (see, e. g., \cite{z06}) and
estimate (\ref{2021-1}), where $l_{1E}=\left\{z=(p_1,p_2,\ldots)^T
\,:\, \|z\|_{1E}\equiv\sum n |p_n| < \infty\right\}$.
\end{proof}

\begin{remark} A number of consequences of this statement can be formulated, for example,
\begin{equation}
\limsup_{t \to \infty}   \|{\bf p}(t)- {\bar{\bf p}}(t)\| \le
\frac{4M\left(M\left|\mathfrak{B-\bar{B}}\right| \mathfrak{f}  +
a\left|\mathfrak{f-\bar{f}}\right|\right)}{ad\left(a -
M\left|\mathfrak{B-\bar{B}}\right|\right)}, \label{3011}
\end{equation}
this follows from
\begin{equation}\|\vp^*-\vp^{**}\| \le 2 \|\vz^*-\vz^{**}\| \le
\frac{4}{d}\|\vz^*-\vz^{**}\|_{1D}.\label{norms2}\end{equation}

The respective perturbation bounds can be formulated for strongly
ergodic (for instance, homogeneous) Markov chains, see
\cite{Zeifman2014s}.

\end{remark}

\begin{remark}
As it was shown in \cite{Zeifman2014s}, the bounds presented in theorem \ref{nonh-pert1} and its corollaries are
sufficiently sharp.
\end{remark}

\section{Convergence rate estimates and perturbation bounds for main classes}

For Markov chains of classes I -- IV an important role is played by the matrix $B^{**}(t)=DB(t)D^{-1}$. To begin with, write out this matrix for each of these classes.

For {\bf class I} this matrix has the form

$B^{**}(t) =  $
\begin{equation}
{=\!\footnotesize \left(\!
\begin{array}{ccccccc}
-\left(\lambda_0\!+\!\mu_1\right)  & \frac{d_1}{d_2}\mu_1
 & 0 & \cdots & 0 &\cdots &\cdots\vspace{1mm}\\
\frac{d_2}{d_1}\lambda_1  & -\!\left(\lambda_1\!+\!\mu_2\right) & \frac{d_2}{d_3}\mu_2 & \cdots & 0 &\cdots &\cdots\vspace{1mm}\\
\ddots & \ddots & \ddots & \ddots & \ddots  &\cdots\vspace{1mm}\\
0 & \cdots & \cdots & \frac{d_{r}}{d_{r-1}}\lambda_{r-1} &
-\!\left(\lambda_{r-1}\!+\!\mu_r\right) & \frac{d_{r}}{d_{r+1}}\mu_{r} &\cdots \vspace{1mm} \\
\cdots & \cdots & \cdots & \cdots & \cdots  & \cdots &\cdots \\
\end{array}
\!\right)} \label{class1-inf}
\end{equation}
\noindent in the case of a countable state space ($S=\infty$);

\smallskip

$B^{**}(t)=$
\begin{equation}
=\!{
\left(\!
\begin{array}{ccccc}
-\left(\lambda_0\!+\!\mu_1\right)  & \frac{d_1}{d_2}\mu_1
 & 0 & \cdots & 0 \vspace{1mm}\\
\frac{d_2}{d_1}\lambda_1  & -\!\left(\lambda_1\!+\!\mu_2\right) & \frac{d_2}{d_3}\mu_2 & \cdots & 0\vspace{1mm} \\
\ddots & \ddots & \ddots & \ddots & \ddots  \vspace{1mm}\\
0 & \cdots & \cdots & \frac{d_{S}}{d_{S-1}}\lambda_{S-1} &
-\!\left(\lambda_{S-1}\!+\!\mu_S\right)
\end{array}
\!\right)} \label{class1-fin}
\end{equation}
\noindent in the case of a finite state space ($S <\infty$).

For {\bf class II} this matrix has the form
\begin{equation}
{
B^{**}(t) =  \left(
\begin{array}{cccccc}
a_{11}  & \frac{d_1}{d_2}\mu_1
 & 0 & \cdots & 0 \vspace{1mm}\\
\frac{d_2}{d_1}a_1  & a_{22} & \frac{d_2}{d_3}\mu_2 & \cdots & 0 \vspace{1mm}\\
\frac{d_3}{d_1}a_2  & \frac{d_3}{d_2}a_1  &a_{33} & \frac{d_3}{d_4}\mu_3 & \cdots &  \vspace{1mm}\\
\ddots & \ddots & \ddots & \ddots & \ddots  \vspace{1mm}\\
\ddots & \ddots & \ddots & \ddots & \ddots  \vspace{1mm}\\
\end{array}
\right)} \label{class2-inf}
\end{equation}
\noindent in the case of a countable state space ($S=\infty$);

\smallskip

$B^{**}(t)=$
\begin{equation}
={
\left(
\begin{array}{ccccc}
a_{11}-a_S  & \frac{d_1}{d_2}\mu_1
 & 0 & \cdots & 0 \vspace{1mm}\\
\frac{d_2}{d_1}\left(a_1-a_S\right)  & a_{22}-a_{S-1} & \frac{d_2}{d_3}\mu_2 & \cdots & 0 \vspace{1mm}\\
\ddots & \ddots & \ddots & \ddots & \ddots  \vspace{1mm}\\
\frac{d_{S}}{d_1}\left(a_{S-1}-a_S\right) & \cdots & \cdots &
\frac{d_{S}}{d_S-1}\left(a_1-a_2\right) & a_{SS}-a_1
\end{array}
\right)} \label{class2-fin}
\end{equation}
\noindent in the case of a finite state space ($S <\infty$).

For {\bf class III} this matrix has the form

\smallskip

$B^{**}(t)=$

\begin{equation}
=\!{\footnotesize \left(\!
\begin{array}{cccccc}
-\left(\lambda_0\! +\!b_1\right)  & \frac{d_1}{d_2}\left(b_1\! -\!
b_2\right)
 &\frac{d_1}{d_3}\left( b_2\! -\! b_3\right) & \cdots & \cdots\vspace{1mm} \\
\frac{d_2}{d_1}\lambda_1 & -\big(\lambda_1\!+\!\sum\limits_{i\le 2}b_i\big) & \frac{d_2}{d_3}\left(b_1\! -\! b_3\right) & \cdots & \cdots\vspace{1mm}\\
\ddots & \ddots & \ddots & \ddots & \ddots\vspace{1mm}  \\
0 & \cdots & \cdots & \frac{d_{r}}{d_{r-1}}\lambda_{r-1} &
-\big(\lambda_{r-1}\!+\!\sum\limits_{i\le r}b_i\big) \cdots\vspace{1mm} \\
\ddots & \ddots & \ddots & \ddots & \ddots  \vspace{1mm}\\
\end{array}
\!\right)}\label{class3-inf}
\end{equation}
\noindent in the case of a countable state space ($S=\infty$);

\smallskip

$B^{**}(t) =$
\begin{equation}
 =\! {\footnotesize  \left(\!
\begin{array}{ccccc}
-\!\left(\lambda_0\! +\!b_1\right)  & \frac{d_1}{d_2}\left(b_1\! -\!
b_2\right)
 & \frac{d_1}{d_3}\left(b_2\! -\! b_3\right) & \cdots & \frac{d_1}{d_S}\left(b_{S-1}\! -\! b_S\right) \vspace{1mm}\\
\frac{d_2}{d_1}\lambda_1  & -\!\big(\lambda_1\!+\!\sum\limits_{i\le 2}b_i\big) & \frac{d_2}{d_3}\left(b_1\! -\! b_3\right) & \cdots & \frac{d_2}{d_S}\left(b_{S-2}\! -\! b_S\right) \vspace{1mm}\\
\ddots & \ddots & \ddots & \ddots & \ddots  \vspace{1mm}\\
0 & \cdots & \cdots & \frac{d_{S}}{d_{S-1}}\lambda_{S-1} &
-\!\big(\lambda_{S-1}\!+\!\sum\limits_{i\le S}b_i\big)
\end{array}
\right)}\label{class3-fin}
\end{equation}
\noindent in the case of a finite state space ($S <\infty$).

Finally, for {\bf class IV} this matrix has the form
\begin{equation}
B^{**} = {
\left(
\begin{array}{cccccc}
a_{11}  & \frac{d_1}{d_2}\left( b_1 - b_2\right)
 & \frac{d_1}{d_3}\left(b_2 - b_3\right) & \cdots & \cdots\vspace{1mm} \\
\frac{d_2}{d_1}a_1  & a_{22} & \frac{d_2}{d_3}\left(b_1 - b_3\right) & \cdots & \cdots \vspace{1mm}\\
\ddots & \ddots & \ddots & \ddots & \ddots  \vspace{1mm}\\
\frac{d_r}{d_1}a_{r-1} & \cdots & \cdots & \frac{d_r}{d_{r-1}}a_1 & a_{rr} & \cdots \vspace{1mm}\\
 \cdots & \cdots & \cdots & \cdots & \cdots & \cdots \\
\end{array}
\right)}\label{class4-inf}
\end{equation}
\noindent in the case of a countable state space ($S=\infty$);

\smallskip

$B^{**}(t)=$
\begin{equation}
=\! {\footnotesize \left(\!
\begin{array}{ccccc}
a_{11}\!-\!a_S  & \frac{d_1}{d_2}\left( b_1\! -\! b_2\right)
 & \frac{d_1}{d_3}\left(b_2\! -\! b_3\right) & \cdots & \frac{d_1}{d_S}\left(b_{S-1}\! -\! b_S\right)\vspace{1mm}\\
\frac{d_2}{d_1}\left(a_1\!-\!a_S\right)  & a_{22}\!-\!a_{S-1} &
\frac{d_2}{d_3}\left(b_1\! -\! b_3\right) & \cdots &
\frac{d_2}{d_S}\left(b_{S-2}\! -\! b_S\right) \vspace{1mm}\\
\ddots & \ddots & \ddots & \ddots & \ddots  \vspace{1mm}\\
\frac{d_{S}}{d_1}\left(a_{S-1}\!-\!a_S\right) & \cdots & \cdots &
\frac{d_{S}}{d_{S-1}}\left(a_1\!-\!a_2\right) & a_{SS}\!-\!a_1
\end{array}\!
\right)}\label{class4-fin}
\end{equation}
\noindent in the case of a finite state space ($S <\infty$).

In the proofs of the following theorems we use the notion of the logarithmic norm of a linear operator function and related estimates of the norm of the Cauchy operator of a linear differential equation. The corresponding results were described in detail in our preceding works, see \cite{gz04,z06,Zeifman2018c}. Here we restrict ourselves only to the necessary minimum.

Recall that the logarithmic norm of an operator function ${B^{**}}(t)$ is defined as the number
$$
\gamma({B^{**}}(t)) = \lim_{h \to
+0}h^{-1}\left(\|I+hB^{**}(t)\|-1\right).
$$
Let $V(t, s)= V(t)V^{-1}(s)$ be the Cauchy operator of the differential equation $$\frac{d{\bf w}}{dt}={B^{**}}(t){\bf w}.$$
Then the estimate
$$\|V(t, s)\| \le e^{\, \int_s^{t} \gamma(B^{**}(u))\, du}$$
holds. Moreover, if for each $t \ge 0$ $B^{**}(t)$ maps $l_1$ into itself, then the logarithmic norm can be calculated by the formula
\begin{equation}
\gamma({B^{**}}(t)) = \sup_{1 \le j \le S}
\left(b_{jj}^{**}(t)+\sum_{i \neq j} |b_{ij}^{**}(t)|\right).
\label{lognorm1}
\end{equation}

Now let
\begin{equation}
\alpha_i\left(t\right)= -\left(b_{jj}^{**}(t)+\sum_{i \neq j}
|b_{ij}^{**}(t)|\right),  \ \alpha\left(t\right)=\inf_{i \ge
1}\alpha_i\left(t\right). \label{posit02}
\end{equation}

Also note that if in classes II--IV the intensities $a_k(t)$ and
$b_k(t)$ do not increase in $k$ for each $t$, then in all the cases
the matrix ${B^{**}}(t)$ is  essentially nonnegative (that is, its
non-diagonal elements are nonnegative), then in (\ref{lognorm1}) and
(\ref{posit02}) the signs of the absolute value can be omitted.

The following statement (\cite[theorem~1]{Zeifman2018c}) is given here for convenience.

\begin{theorem}{}
Let for some sequence $\{d_i, \ i \ge 1\}$ of positive numbers the conditions $d_1=1$, $d=\inf_{i \ge 1} d_i > 0$ and
\begin{equation}
\int_0^{\infty} \alpha(t)\, dt = + \infty \label{ord-erg1}
\end{equation}
hold. Then the Markov chain $X(t)$ is weakly ergodic and for any initial condition $s \ge 0$, ${\bf w}(s)$ and for all $t \ge s$ there holds the estimate
\begin{eqnarray}
\|{\bf w}\left(t\right)\| \le e^{-\int_s^t
{\alpha\left(u\right)du}}\|{\bf w}(s)\|. \label{t001}
\end{eqnarray}
\end{theorem}

Now let instead of (\ref{ord-erg1}), for all $0 \le s \le t$ a stronger condition
\begin{equation}
e^{-\int_s^{t} \alpha(\tau)\, d\tau} \le M^*e^{-a^*(t-s)}
\label{ord-erg2}
\end{equation}
holds.

\begin{theorem}{}
Let under the conditions of theorem 3 inequality (\ref{ord-erg2}) holds. Then the Markov chain $X(t)$ is $1D$-exponentially weakly ergodic and for all $t\ge s \ge 0$ and ${\bf p^{*}}(s) \in l_{1D}$, ${\bf p^{**}}(s) \in l_{1D}$ there holds the inequality (\ref{d-erg2}) with $M=M^*$ and $a=a^*$.
\end{theorem}

\begin{remark} In the case of a homogeneous Markov chain or if all intensities are periodic with one and the same period, conditions (\ref{ord-erg1}) and (\ref{ord-erg2}) are equivalent.
\end{remark}

\begin{theorem}{}
Let the conditions of theorem 4 hold. Then the Markov chain $X(t)$ is $1D$-exponentially weakly ergodic, under perturbations small enough (\ref{3001'}) the perturbed chain $\bar{X}(t)$ is also $1D$-exponentially weakly ergodic and perturbation bound (\ref{3002}) in the $1D$-norm holds. If, moreover, $W=\inf_{i \ge 1} \frac {d_i}{i} > 0$, then both chains $X(t)$ and $\bar{X}(t)$ have limit expectations and estimate (\ref{3003}) holds for the perturbation of the mathematical expectation.
\end{theorem}

To obtain perturbation bounds in the natural norm it suffices to use inequality (\ref{norms2}) mentioned above.

\begin{corollary}
Under the conditions of theorem 5 the following perturbation bound in the natural $l_1$- (total variation) norm holds:
\begin{eqnarray}
\limsup_{t \to \infty}   \|{\bf p}(t)- {\bar{\bf p}}(t)\| \le
\frac{4M\left(M\left|\mathfrak{B-\bar{B}}\right| \mathfrak{f}  +
a\left|\mathfrak{f-\bar{f}}\right|\right)}{ad\left(a -
M\left|\mathfrak{B-\bar{B}}\right|\right)}.
\label{3002}\end{eqnarray}
\end{corollary}

Note that it is convenient to use the results formulated above for the construction of perturbation bounds for Markov chains of the first four classes, see, e. g., \cite{z98,Zeifman2014s,Zeifman2014q,Zeifman2018c}.

For chains of the fifth class, as a rule, it is convenient to use the approach based on uniform bounds as will be shown below. These models were considered, e. g., in \cite{z12,Zeifman2017sm,Zeifman2017ecms}.

Let
\begin{equation}
\beta_*\left(t\right) = \inf_k a_{0k}(t).\label{0001}
\end{equation}

\begin{theorem} Let the intensities of catastrophes be essential, that is
\begin{equation}
\label{bet*} \int_0^\infty \beta_*\left(t\right)\, dt = +\infty.
\end{equation}
Then the chain $X\left(t\right)$ is weakly ergodic in the uniform operator topology and for any initial conditions ${{\bf
p}}^{*}\left(0\right), {{\bf p}}^{**}\left(0\right)$ and any $0 \le s \le t$ there holds the following convergence rate estimate:
\begin{eqnarray}
\left\|{{\bf p}}^{*}\left(t\right)-{{\bf
p}}^{**}\left(t\right)\right\|  \le 2 e^{-\int\limits_s^t
\beta_*\left(\tau\right)\, d\tau}. \label{0011}\end{eqnarray}
\end{theorem}

To prove this theorem we will use the same technique as in \cite{z12}. Rewrite the direct Kolmogorov system (\ref{ur01}) in the form
\begin{equation}
\frac{d {\bf p} }{dt}=A^*\left( t\right) {{\bf p} }  +{\bf g}
\left(t\right), \quad t\ge 0. \label{eq0112''}
\end{equation}
\noindent Here ${\bf
g}\left(t\right)=\left(\beta_*\left(t\right),0,0, \dots\right)^T$,
$A^*\left(t\right)=\left (
a_{ij}^*\left(t\right)\right)_{i,j=0}^{\infty}$, and
\begin{equation}
a_{ij}^*\left(t\right) = \begin{cases}
a_{0j}\left(t\right) - \beta_*\left(t\right), & \mbox { if }  i= 0, \vspace{1mm}\cr
a_{ij}\left(t\right), & \mbox { otherwise }.
\end{cases}
\label{01101}
\end{equation}
The solution to this equation can be written as
\begin{equation}
{\bf p}\left(t\right)= U^*\left(t,0\right){\bf
p}\left(0\right)+\int_0^t{}U^*\left(t,\tau\right){\bf
g}\left(\tau\right)\,d\tau, \label{0appr218}
\end{equation}
\noindent where $U^*\left(t,s\right)$ is the Cauchy operator of the differential equation
\begin{equation}
\frac{d {\bf z} }{dt}=A^*\left( t\right) {{\bf z} }.
\label{0eq112'''}
\end{equation}
Note that the matrix $A^*\left( t\right)$ is essentially nonnegative
for all $t \ge 0$. Its logarithmic norm is equal to
\begin{equation}
\gamma(A^*(t)) = \sup_i \left(a_{ii}^*\left(t\right) + \sum_{j\neq
i} a_{ji}^*\left(t\right)\right) = -\beta_*\left(t\right),
\label{0cat03'}
\end{equation}
\noindent and hence,
\begin{eqnarray}
\left\|{{\bf p} }^{*}\left(t\right)-{{\bf p}
}^{**}\left(t\right)\right\|\le e^{-\int\limits_s^t
\beta_*\left(\tau\right)\, d\tau} \left\|{{\bf p}
}^{*}\left(s\right)-{{\bf p} }^{**}\left(s\right)\right\| \le
2e^{-\int\limits_s^t \beta_*\left(\tau\right)\, d\tau}.
\label{bound3'}
\end{eqnarray}

\begin{theorem} Let, instead of (\ref{bet*}), the stronger condition
\begin{equation}
\label{bet**} e^{-\int_s^{t} \beta_*(\tau)\, d\tau} \le
c^*e^{-b^*(t-s)}
\end{equation}
hold. Then the chain $X\left(t\right)$ is weakly exponentially ergodic in the uniform operator topology, and if the perturbations are small enough, that is, $\|\hat{A}(t)\| \le \varepsilon$ for almost all $t \ge 0$, then the perturbed chain $\bar{X}(t)$ is also exponentially weakly ergodic and the perturbation bound (\ref{207''}) holds with $c=c^*$, $b=b^*$.
\end{theorem}

\section{Examples}

First of all note that many examples of perturbation bounds for queueing systems were considered in
\cite{z94,z98,Zeifman2010icumt,Zeifman2011value,z12,Zeifman2014q,Zeifman2014s}.

Here, to compare both approaches, we will mostly deal with the
queueing system $M_t|M_t|N|N$ with losses and $1$-periodic
intensities. In the preceding papers on this model, other problems
were considered. For example, in \cite{Doorn2010} the asymptotics of
the rate of convergence to the stationary mode as $N \to \infty$,
was studied, whereas the paper \cite{Doorn2011} dealt with the
asymptotics of the convergence parameter under various limit
relations between the intensities and the dimensionality of the
model. In \cite{Zeifman2010icumt,Zeifman2011value} perturbation
bounds were considered under additional assumptions.

Let $N \geq 1$ be the number of servers in the system.  Assume that
the customers arrival intensity $\la(t)$ and the service intensity
of a server $\mu(t)$ are $1$-periodic nonnegative functions
integrable on the interval $[0,1]$. Then the number of customers in
the system (queue length) $X(t)$ is a finite Markov chain of class
I, that is, a BDP with the intensities $\la_{k-1}(t)= \la(t)$,
$\mu_k(t) = k \mu(t)$ for $k=1,\dots,N$.

It should be especially  noted that the process $X(t)$ is weakly
ergodic (obviously, exponentially and uniformly ergodic, since the
intensities are periodic and the state space is finite), if and only
if
\begin{equation}
\int_0^1 \left(\la (t)+ \mu (t)\right)\, dt >0,\label{example000}
\end{equation}
\noindent see, e. g., \cite{Zeifman2009ait}.

For definiteness, assume that $\int_0^1 \mu (t)\, dt>0.$

Apply the approach described in theorems 3 and 4.

Let all $d_k=1$. Then
\begin{equation}
B^{**}(t)={
\left(
\begin{array}{ccccc}
-\left(\lambda+\mu\right)  & \mu
 & 0 & \cdots & 0 \\
\lambda  & -\left(\lambda+2\mu\right) & 2\mu & \cdots & 0 \\
\ddots & \ddots & \ddots & \ddots & \ddots  \\
0 & \cdots & \cdots & \lambda & -\left(\lambda+N\mu\right)
\end{array}
\right),} \label{example001}
\end{equation}
\noindent and in (\ref{posit02}) we have $\alpha_i\left(t\right)=\mu(t)$ for all $i$, hence, $\alpha\left(t\right)=\mu(t)$.

Therefore, theorem 3 yields the estimate
\begin{equation}
\|{\bf p^{*}}(t) - {\bf p^{**}}(t)\|_{1D} \le  e^{-\int_s^t  \mu
(\tau)\, d\tau} \|{\bf p^{*}}(s) - {\bf p^{**}}(s)\|_{1D}.
\label{example002} \end{equation}
To find the constants in the estimates, let $\mu^*=\int_0^1 \mu(\tau)\, d\tau $ and consider
\begin{equation}
\int_0^t  \mu (\tau)\, d\tau = \mu^* t + \int_0^{\{t\}} \left( \mu
(\tau)- \mu^*\right)\, d\tau.  \label{example003}
\end{equation}
Find the bound for the second summand in (\ref{example003}). Assuming $u=\{t\}$, we obtain
\begin{equation}
\left| \int_0^{u} \left( \mu (\tau)- \mu^*\right)\, d\tau\right| \le
K^* = \sup_{u \in [0,1]}\int_0^{u} \left( \mu (\tau)- \mu^*\right)\,
d\tau. \label{example004}
\end{equation}
Then
\begin{equation}
e^{-\int_s^t  \mu (\tau)\, d\tau} \le e^{K^*}
e^{-\mu^*\left(t-s\right)}. \label{example005}
\end{equation}
Therefore, for the queueing system $M_t|M_t|N|N$ the conditions of theorem 5 and corollary 5 hold with
\begin{equation}d=1,\quad M=M^*=e^{K^*}, \quad a=a^* = \mu^*, \quad W=\frac{1}{N}. \label{example006}
\end{equation}
These statements imply the following perturbation bounds:
\begin{eqnarray}
\limsup_{t \to \infty}   \|{\bf p}(t)- {\bar{\bf p}}(t)\| \le
\frac{4e^{K^*}\left(e^{K^*}\left|\mathfrak{B-\bar{B}}\right|
\mathfrak{f} +
\mu^*\left|\mathfrak{f-\bar{f}}\right|\right)}{\mu^*\left(\mu^* -
e^{K^*}\left|\mathfrak{B-\bar{B}}\right|\right)}, \label{example021}
\end{eqnarray}
\noindent for the vector od=f state probabilities and
\begin{eqnarray}
\limsup_{t \to \infty}  |\phi(t) - \bar{\phi}(t)|\le
\frac{Ne^{K^*}\left(e^{K^*}\left|\mathfrak{B-\bar{B}}\right|
\mathfrak{f} +
\mu^*\left|\mathfrak{f-\bar{f}}\right|\right)}{\mu^*\left(\mu^* -
e^{K^*}\left|\mathfrak{B-\bar{B}}\right|\right)}, \label{example022}
\end{eqnarray}
\noindent for limit expectations.

Moreover, for these bounds to be consistent, the additional in=formation is required concerning the form of the perturbed intensity matrix. The simplest bounds can be obtained, if it is assumed that the perturbed Markov chain is also a BDP with the same state space and birth and death intensities $\la_{k-1}(t)$ and $\mu_k(t)$ respectively. Then if the birth and death intensities themselves do not exceed $\varepsilon$ for almost all $t \ge 0$, then $|\mathfrak{f-\bar{f}}| \le \varepsilon$ and $|\mathfrak{B-\bar{B}}| \le 5 \varepsilon$, so that the bounds
(\ref{example021}) and (\ref{example022}) have the form
\begin{eqnarray}
\limsup_{t \to \infty}   \|{\bf p}(t)- {\bar{\bf p}}(t)\| \le
\frac{4e^{K^*}\left(5L e^{K^*}  +
\mu^*\right)\varepsilon}{\mu^*\left(\mu^* - 5\varepsilon
e^{K^*}\right)}, \label{example023}
\end{eqnarray}
\noindent for the vectors of state probabilities and
\begin{eqnarray}
\limsup_{t \to \infty}  |\phi(t) - \bar{\phi}(t)|\le
\frac{4Ne^{K^*}\left(5L e^{K^*}  +
\mu^*\right)\varepsilon}{\mu^*\left(\mu^* - 5\varepsilon
e^{K^*}\right)}, \label{example024}
\end{eqnarray}
\noindent for the limit expectations.

On the other hand, theorem 7 can be applied as well. To construct the bounds for the corresponding parameters, use (\ref{norms2}) and the fact that
$\|D\|_1=N$. Then theorem 7 is valid for the queueing system $M_t|M_t|N|N$ with the following values of the parameters:
\begin{equation}
c=c^*=4Ne^{K^*}, \quad b=b^* = \mu^*. \label{example007}
\end{equation}
According to this theorem we obtain the estimate
\begin{equation}
\limsup_{t \to \infty}   \|{\bf p}(t)- {\bar{\bf p}}(t)\| \le
\frac{\left(1+ K^*+ \log (2N)\right)\varepsilon}{\mu^*}.
\label{example011}
\end{equation}
Moreover, the Markov chains $X(t)$ and $\bar{X}(t)$ have limit expectations and
\begin{equation}
|\phi (t) - \bar{\phi} (t)| \le  \frac{N\left(1+K^*+ \log
(2N)\right)\varepsilon}{\mu^*} . \label{example012}
\end{equation}

It is worth noting that for estimates (\ref{example011}) and (\ref{example012}) to hold, only the condition of the smallness of perturbations is required and {\it no} additional information concerning the structure of the intensity matrix is required.

Thus, in the example with the finite state space uns=der consideration, uniform bounds turn out to be more exact.

\smallskip

Now consider a more special example. Let $N=299$, $\la(t)=200(1+\sin 2\pi \omega t)$, $\mu(t) =1$.

On Fig. 1--5 there are the plots of the expected  number of
customers in the system for some of most probable states with
$\omega=1$; on Fig. 6--7 there are the plots of of the expected
number of customers with $\omega=0.5$.

On the other hand, as it has already been noted, for the Markov chains of classes I--IV with countable state space no uniform bounds could be constructed.

Consider the construction of bounds on the example of a rather simple model, which, however, does not belong to the most well-studied class I (that is, which is not a BDP).

Let a queueing system be given in which the  customers can appear
separately or in pairs with the corresponding intensities
$a_1(t)=\la(t)$ and $a_2(t)=0.5\la(t)$, but are served one by one on
one of two servers with constant intensities
$\mu_k(t)=\min(k,2)\mu$, where $\la(t)$ is a $1$-periodic function
integrable on the interval $[0,1]$. Then the number of customers in
this system belongs to class II and the corresponding matrix
$B^{**}(t)$ has the form
\begin{equation}
{
B^{**}(t) =  \left(
\begin{array}{cccccc}
a_{11}  & \frac{d_1}{d_2}\mu
 & 0 & \cdots & 0 \vspace{1mm}\\
\frac{d_2}{d_1}\la  & a_{22} & \frac{d_2}{d_3}2\mu & \cdots & 0 \vspace{1mm}\\
\frac{d_3}{d_1}0.5\la  & \frac{d_3}{d_2}\la  &a_{33} & \frac{d_3}{d_4}2\mu & \cdots &  \\
0 & \ddots & \ddots & \ddots & \ddots  \\
\ddots & \ddots & \ddots & \ddots & \ddots  \\
\end{array}
\right)}, \label{example031}
\end{equation}
\noindent where $a_{11}(t)= -\left(1.5\la(t)+\mu\right)$, $a_{kk}(t)=
-\left(1.5\la(t)+2\mu\right)$, if $k \ge 2$.
This matrix is essentially nonnegative, so that in the expression for the logarithmic norm the signs of the absolute value can be omitted. Let $d_1=1$, $d_{k+1}=\delta d_k$ and choose $\delta>1$. For this purpose consider the expressions from (\ref{posit02}). We have
$$\alpha_1(t)=\mu-\la(t)\left(0.5\delta ^2+\delta - 1.5\right),$$
$$\alpha_2(t)=\mu\left(2- \delta^{-1}\right)-\la(t)\left(0.5\delta
^2+\delta - 1.5\right),$$
$$\alpha_k(t)=2\mu\left(1- \delta^{-1}\right)-\la(t)\left(0.5\delta
^2+\delta - 1.5\right),\ k \ge 3.$$
Then for $\delta \le 2$ we obtain
\begin{eqnarray}\alpha\left(t\right)=\inf_{i \ge
1}\alpha_i\left(t\right)=2\mu\left(1-
\delta^{-1}\right)-\la(t)\left(0.5\delta ^2+\delta - 1.5\right)=
\nonumber\\ = \left(\delta - 1\right)\left(\frac{2\mu}{
\delta}-0.5\la(t)\left(\delta +3 \right)\right), \label{example032}
\end{eqnarray}
\noindent and the condition
\begin{equation}
\alpha^* = \int_0^1 \left(\delta - 1\right)\left(\frac{2\mu}{
\delta}-0.5\la(t)\left(\delta +3 \right)\right)\, dt = \frac{\delta
- 1}{2}\left(\frac{4\mu}{ \delta}-\la^*\left(\delta +3
\right)\right)
>0 , \label{example033}
\end{equation}
\noindent will a fortiori hold, if $\mu > \la^*$ with a corresponding choice of $\delta \in (1,2]$.

The further reasoning is almost the same as in the preceding example: instead of (\ref{example005}) we obtain
\begin{equation}
e^{-\int_s^t  \alpha (\tau)\, d\tau} \le e^{K^*}
e^{-\alpha^*\left(t-s\right)}, \label{example034}
\end{equation}
\noindent where now
\begin{equation}
K^* = \sup_{u \in [0,1]}\int_0^{u} \left( \alpha (\tau)-
\alpha^*\right)\, d\tau. \label{example035}
\end{equation}
Hence, the conditions of theorem 5 and corollary 5 for the number of
customers in the system under consideration hold for
\begin{equation}d=1,\quad M=M^*=e^{K^*}, \quad a=a^* = \alpha^*, \quad W=\inf_{k \ge 1}\frac{\delta^{k-1}}{k}. \label{example036}
\end{equation}

To  construct meaningful perturbation bounds, it is necessarily
required to have additional information concerning the form of the
perturbed intensity matrix. So, example 1 in Section 2 shows that if
a possibility of the arrival of an arbitrary number of customers
(`mass arrival' in the terminology of \cite{Zhang2015}) to an empty
queue is assumed, then an arbitrarily small (in the uniform norm)
perturbation of the intensity matrix can `spoil' all the
characteristics of the process. For example, satisfactory bounds can
be constructed, if we know that the intensity matrix of the
perturbed system has the same form, that is, the customers can
appear either separately or in pairs and are served one by one.
Then, if the perturbations of the intensities themselves do not
exceed $\varepsilon$ for almost all $t \ge 0$, then
$|\mathfrak{f-\bar{f}}| \le 5\varepsilon$ and
$|\mathfrak{B-\bar{B}}| \le 5 \varepsilon$, so that instead of
(\ref{example021}) and (\ref{example022}) we obtain
\begin{eqnarray}
\limsup_{t \to \infty}   \|{\bf p}(t)- {\bar{\bf p}}(t)\| \le
\frac{20e^{K^*}\varepsilon\left(L e^{K^*}  +
\alpha^*\right)}{\alpha^*\left(\alpha^* - 20\varepsilon
e^{K^*}\right)}, \label{example039}
\end{eqnarray}
\noindent for the vectors of state probabilities and
\begin{eqnarray}
\limsup_{t \to \infty}  |\phi(t) - \bar{\phi}(t)|\le
\frac{20e^{K^*}\varepsilon\left(L e^{K^*}  +
\alpha^*\right)}{\alpha^*W\left(\alpha^* - 20\varepsilon
e^{K^*}\right)}, \label{example040}
\end{eqnarray}
\noindent for the limit expectations.

The particular example: let $\la(t)=1+\sin 2\pi t$, $\mu(t) =3$. Choose $\delta=2$. Then we have
\begin{equation}\alpha(t)= \mu-2.5\la(t), \quad  \alpha^*=0.5, \quad W=1. \label{example041}
\end{equation}

Further we follow the method that was described in
\cite{Zeifman2014i,Zeifman2017tpa} in detail. Namely, we choose the
dimensionality of the truncated process (300 in our case), the
interval on which the desired accuracy is achieved ($[0,100]$) in
the example under consideration) and the limit interval itself (here
it is $[100,101]$).

Fig. 8--13 expose the plots of the expected number of customers in
the system and some most probable states.

\section*{Acknowledgement}

Sections 1--3 were written by Korolev and Zeifman under the support of the Russian Science Foundation, project 18-11-00155, Sections 4 and 5 were written by Satin and Zeifman under the support of the Russian Science Foundation, project 19-11-00020.

\renewcommand{\refname}{References}

\renewcommand{\figurename}{Fig.}

\begin{figure}
\begin{center}
\includegraphics[width=14cm]{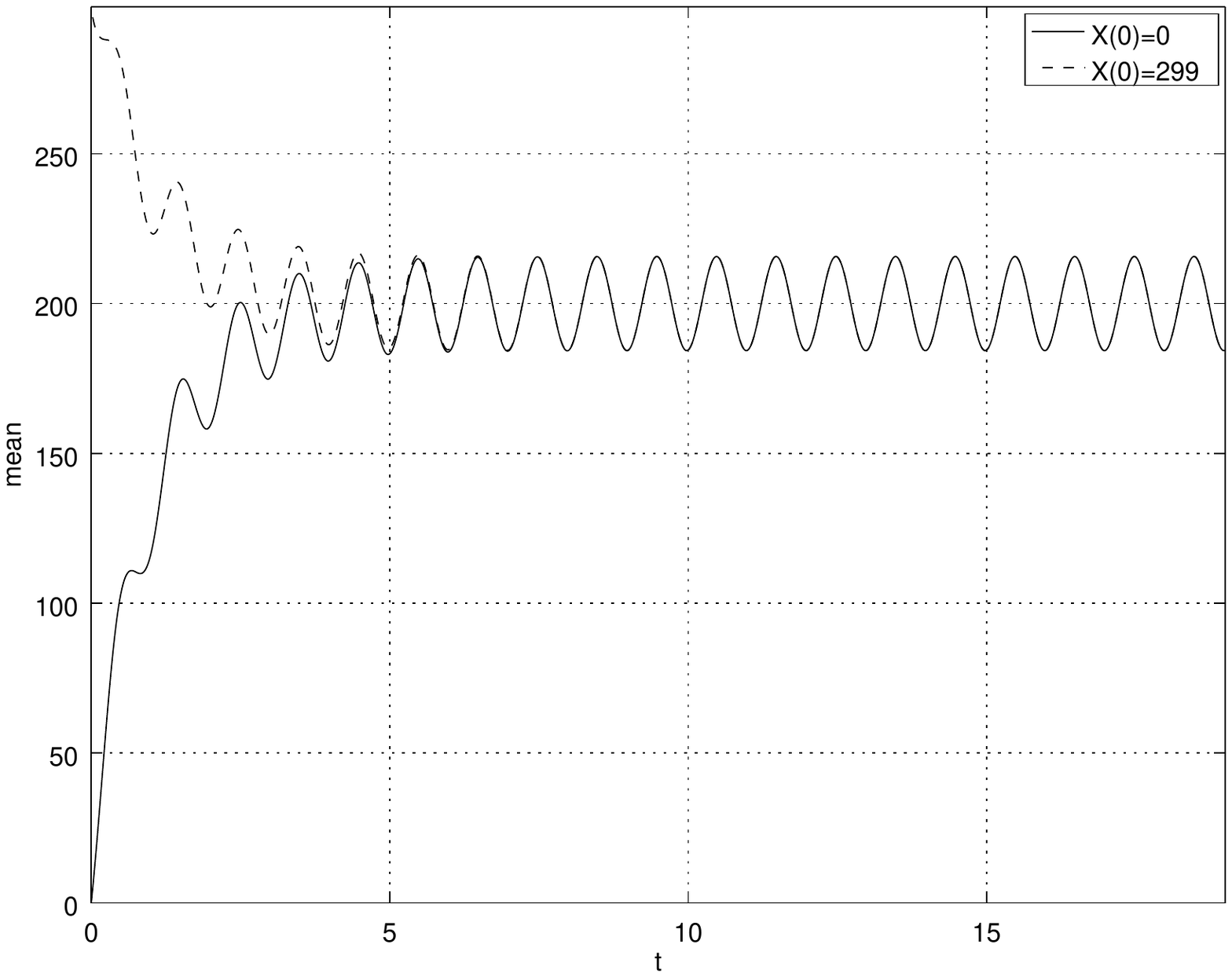}
\end{center}
\vspace{-5cm}\caption{Example 1. The mean $E(t,0)$ and $E(t,N)$ for
the original process $t\in[0,19]$, $\omega=1$. }
\end{figure}

\begin{figure}
\begin{center}
\includegraphics[width=14cm]{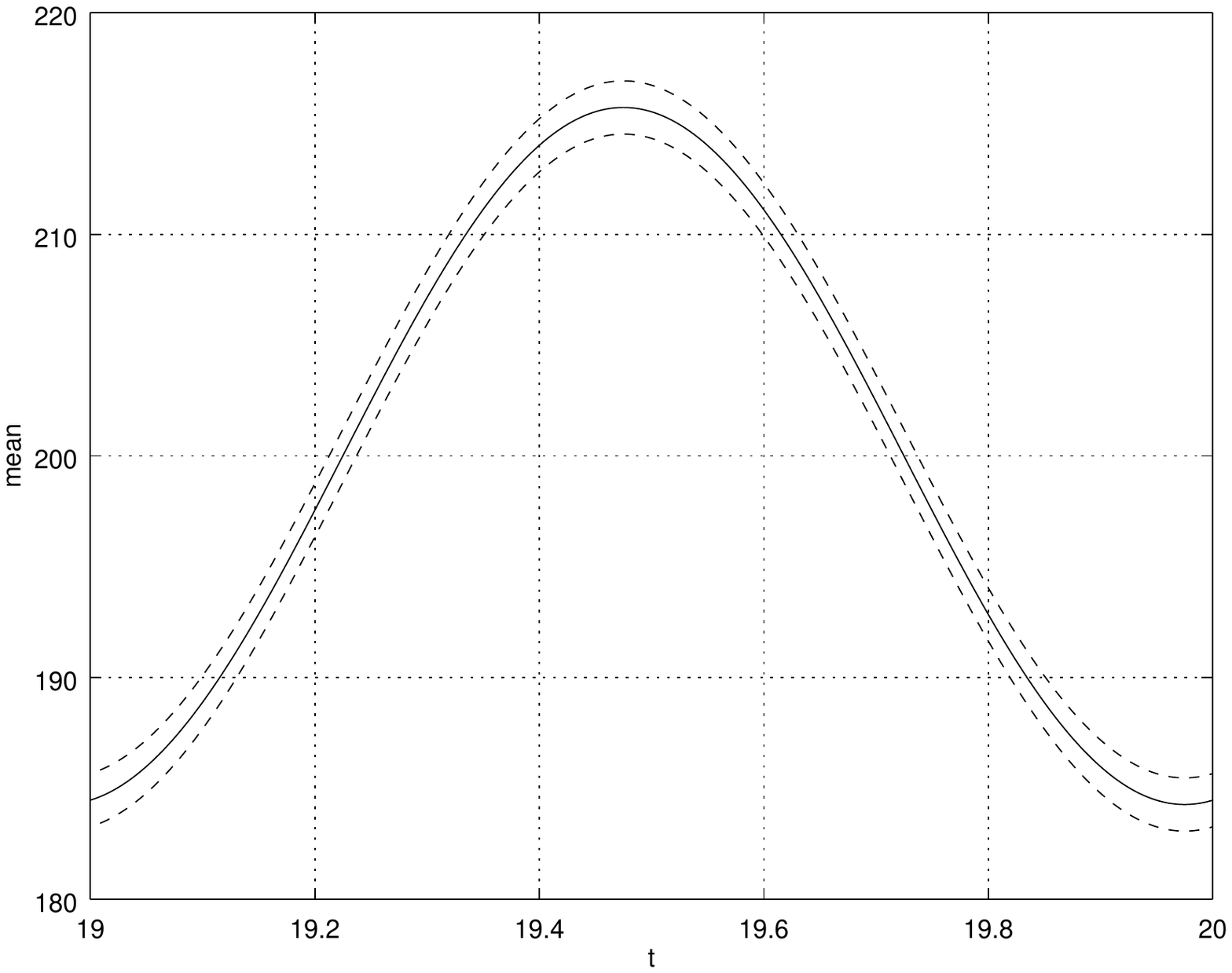}
\end{center}
\vspace{-5cm}\caption{Example 1. The perturbation bounds for the
limit expectation $E(t,0)$, $t\in[19,20]$, $\omega=1$. }
\end{figure}

\begin{figure}
\begin{center}
\includegraphics[width=14cm]{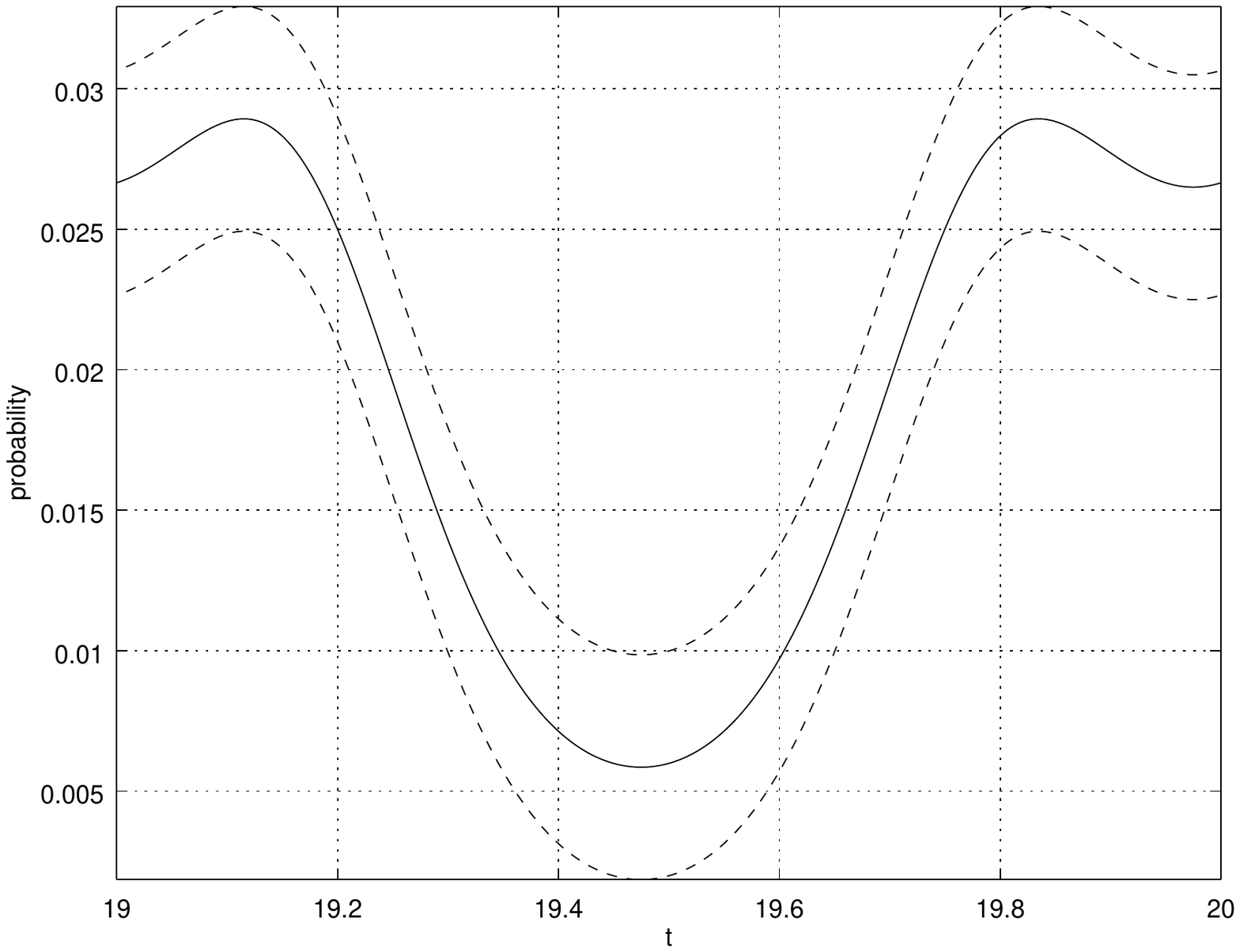}
\end{center}
\vspace{-5cm}\caption{Example 1. The perturbation bounds for the
'limit' probability $\Pr(X(t)=190)$, $t\in[19,20]$, $\omega=1$. }
\end{figure}

\begin{figure}
\begin{center}
\includegraphics[width=14cm]{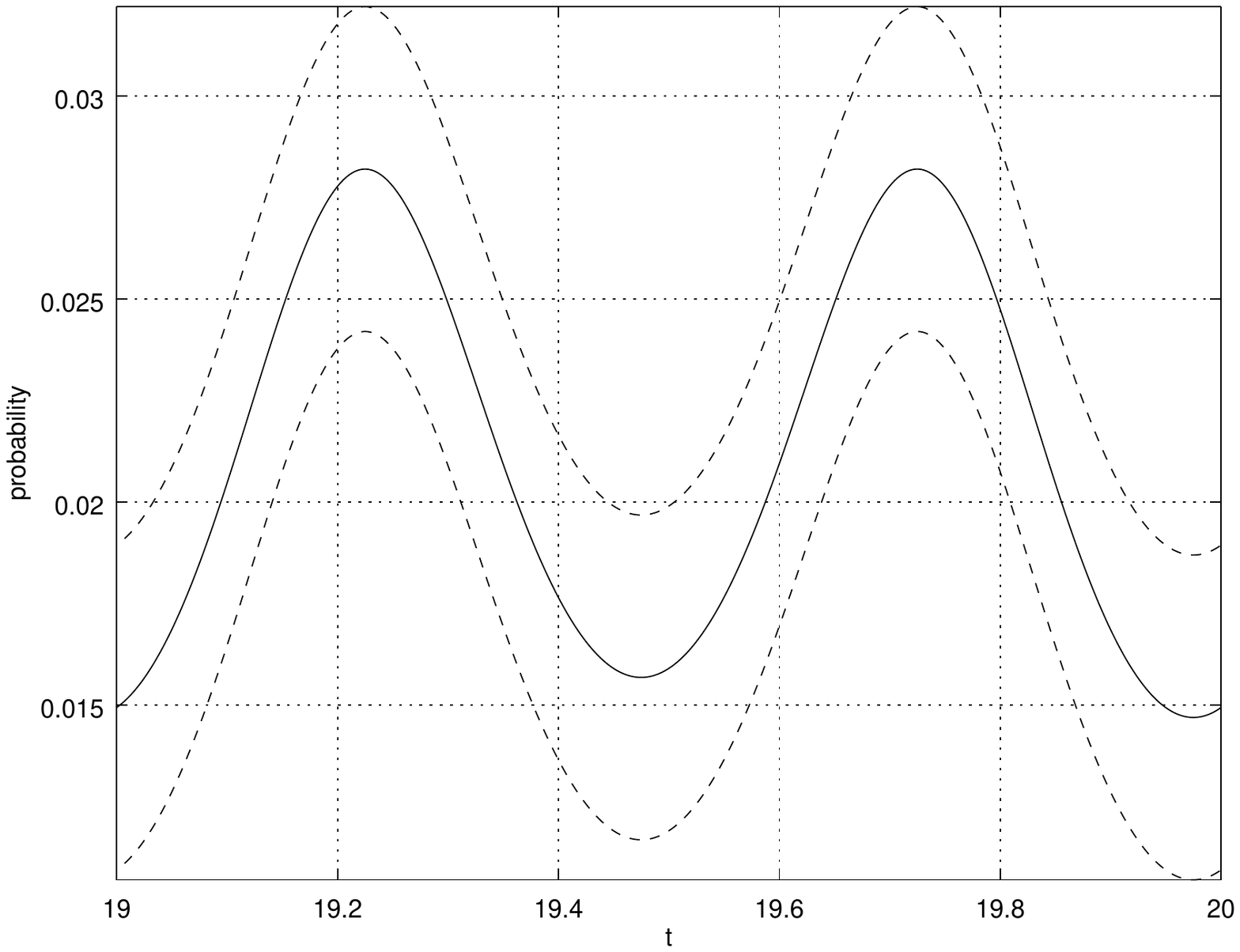}
\end{center}
\vspace{-5cm}\caption{Example 1. The perturbation bounds for the
'limit' probability $\Pr(X(t)=200)$, $t\in[19,20]$, $\omega=1$. }
\end{figure}

\begin{figure}
\begin{center}
\includegraphics[width=14cm]{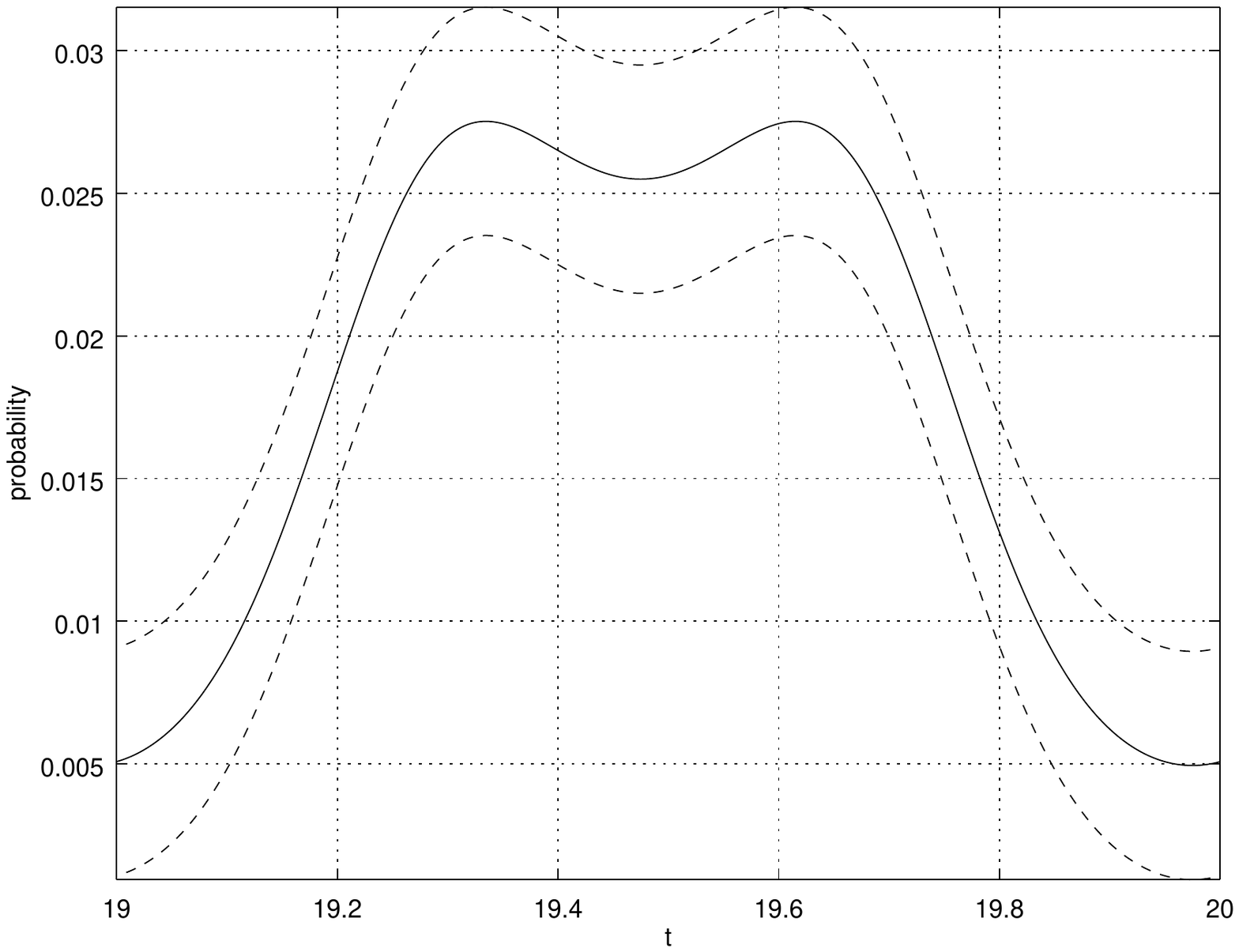}
\end{center}
\vspace{-5cm}\caption{Example 1. The perturbation bounds for the
'limit' probability $\Pr(X(t)=210)$, $t\in[19,20]$, $\omega=1$. }
\end{figure}

\begin{figure}
\begin{center}
\includegraphics[width=14cm]{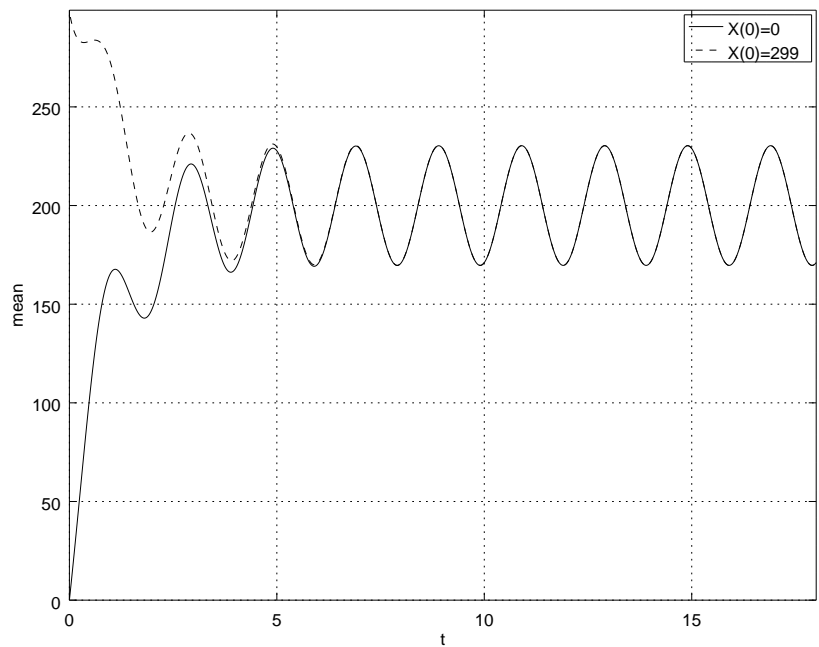}
\end{center}
\vspace{-5cm}\caption{Example 1. The expectations $E(t,0)$ and
$E(t,N)$ for the original process $t\in[0,18]$, $\omega=0.5$. }
\end{figure}

\begin{figure}
\begin{center}
\includegraphics[width=14cm]{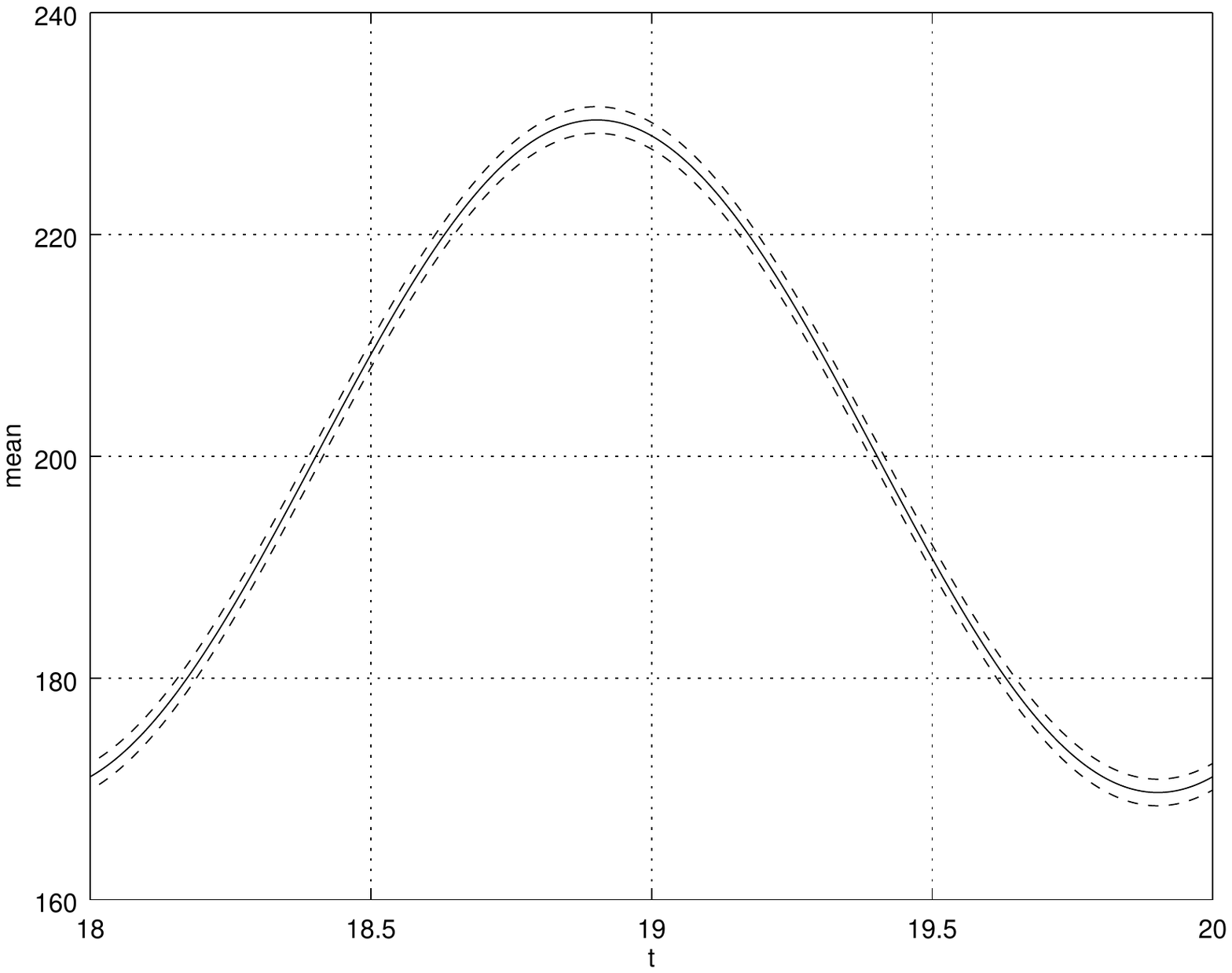}
\end{center}
\vspace{-5cm}\caption{Example 1. The perturbation bounds for the
limit expectation $E(t,0)$, $t\in[18,20]$, $\omega=0.5$. }
\end{figure}

\begin{figure}
\begin{center}
\includegraphics[width=14cm]{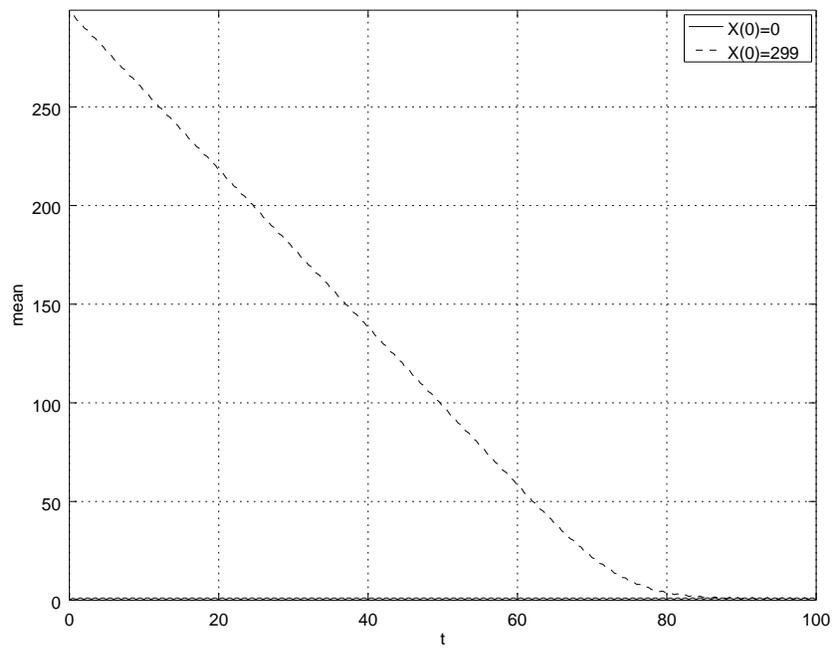}
\end{center}
\vspace{-5cm}\caption{Example 2. The expectations $E(t,0)$ and
$E(t,299)$ for the original process $t\in[0,100]$. }
\end{figure}

\begin{figure}
\begin{center}
\includegraphics[width=14cm]{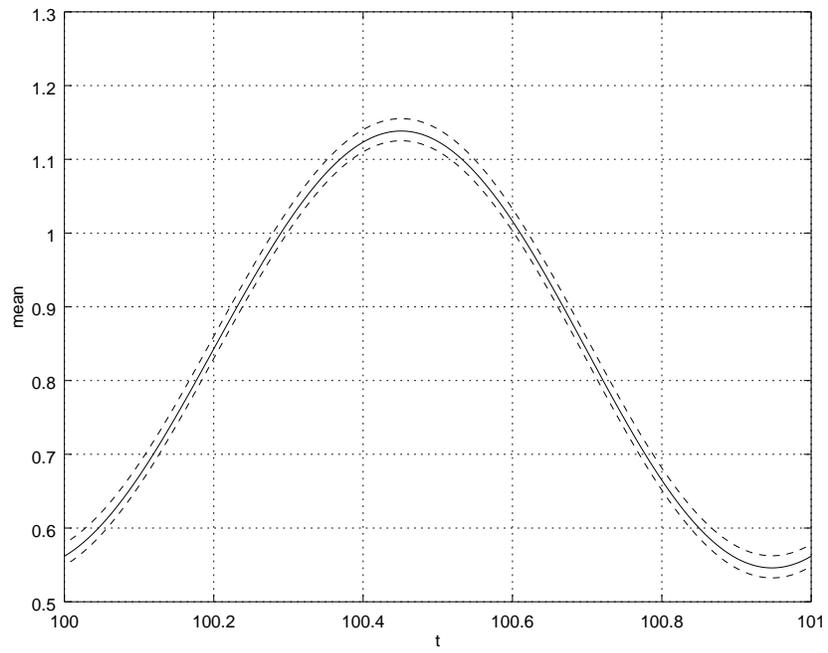}
\end{center}
\vspace{-5cm}\caption{Example 2.  The perturbation bounds for the
limit expectation $E(t,0)$, $t\in[100,101]$. }
\end{figure}

\begin{figure}
\begin{center}
\includegraphics[width=14cm]{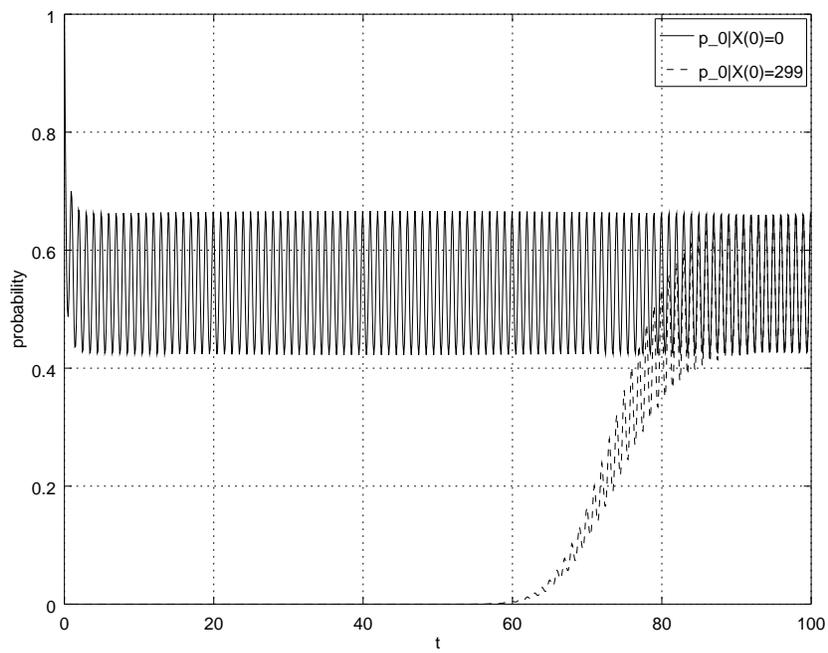}
\end{center}
\vspace{-5cm}\caption{Example 2. The probabilities of the empty
queue for $X(0)=0$ and $X(0)=299$ for the original process
$t\in[0,100]$. }
\end{figure}

\begin{figure}
\begin{center}
\includegraphics[width=14cm]{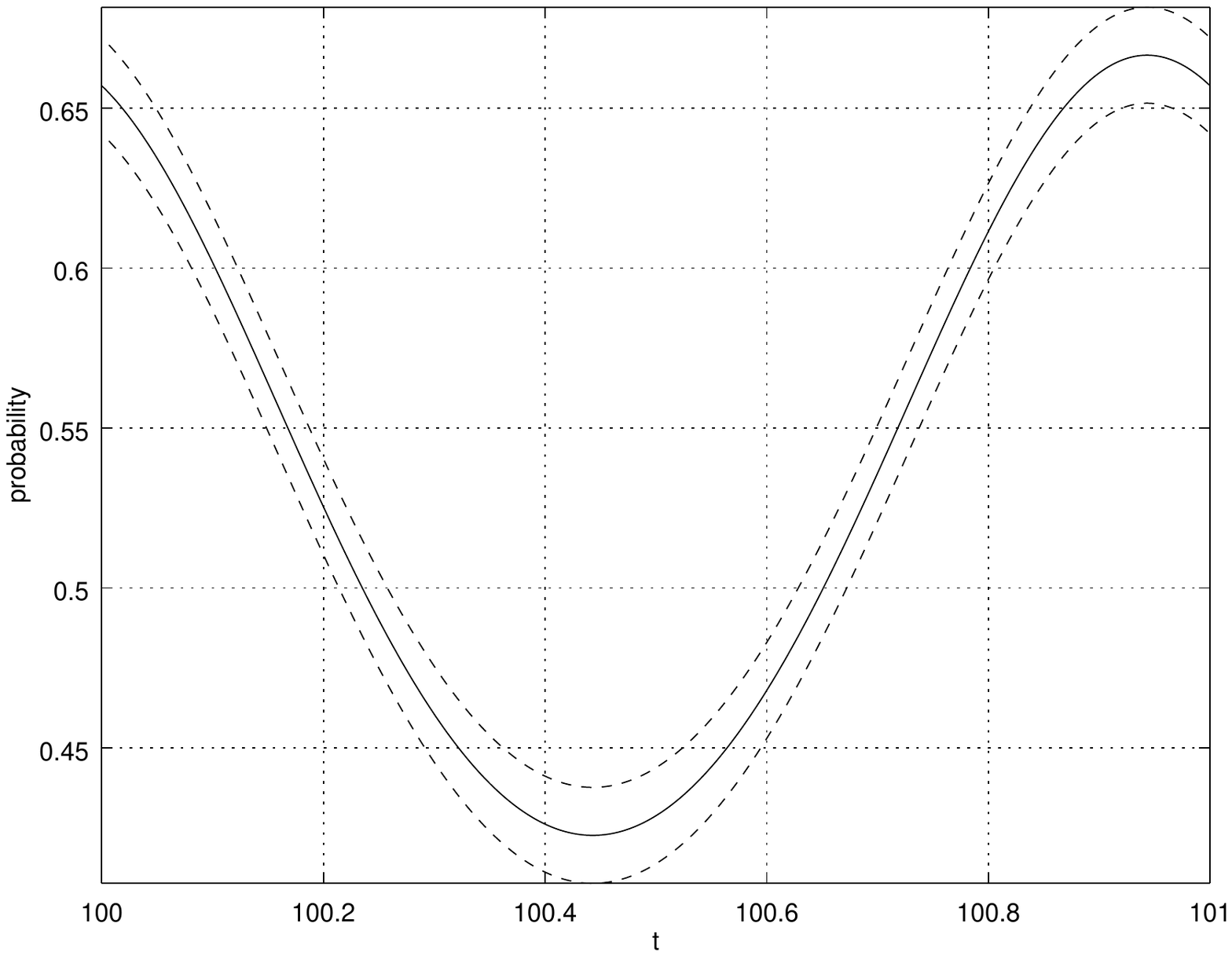}
\end{center}
\vspace{-5cm}\caption{Example 2. The perturbation bounds for the
'limit' probability $\Pr(X(t)=0)$, $t\in[100,101]$. }
\end{figure}

\begin{figure}
\begin{center}
\includegraphics[width=14cm]{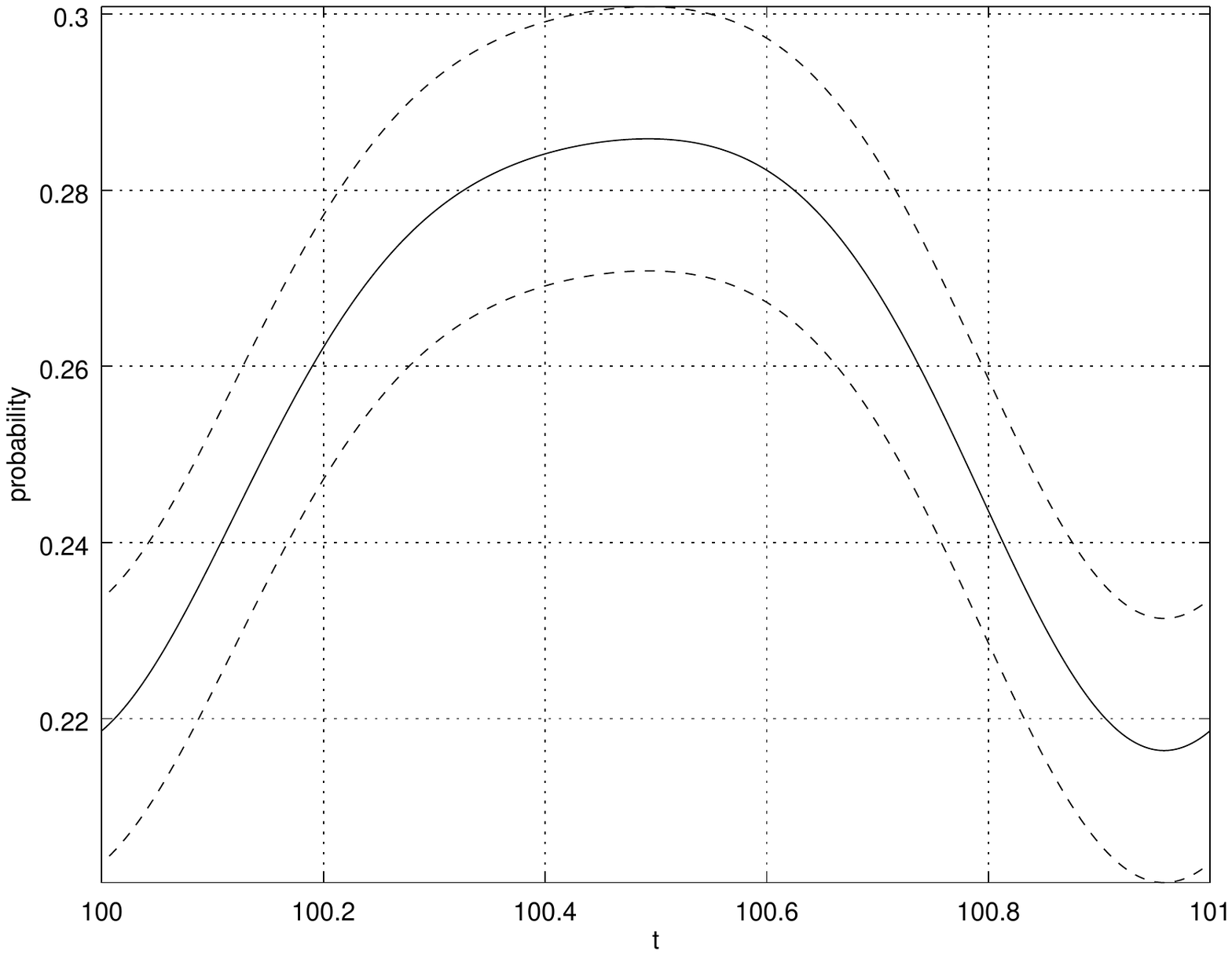}
\end{center}
\vspace{-5cm}\caption{Example 2. The perturbation bounds for the
'limit' probability $\Pr(X(t)=1)$, $t\in[100,101]$. }
\end{figure}

\begin{figure}
\begin{center}
\includegraphics[width=14cm]{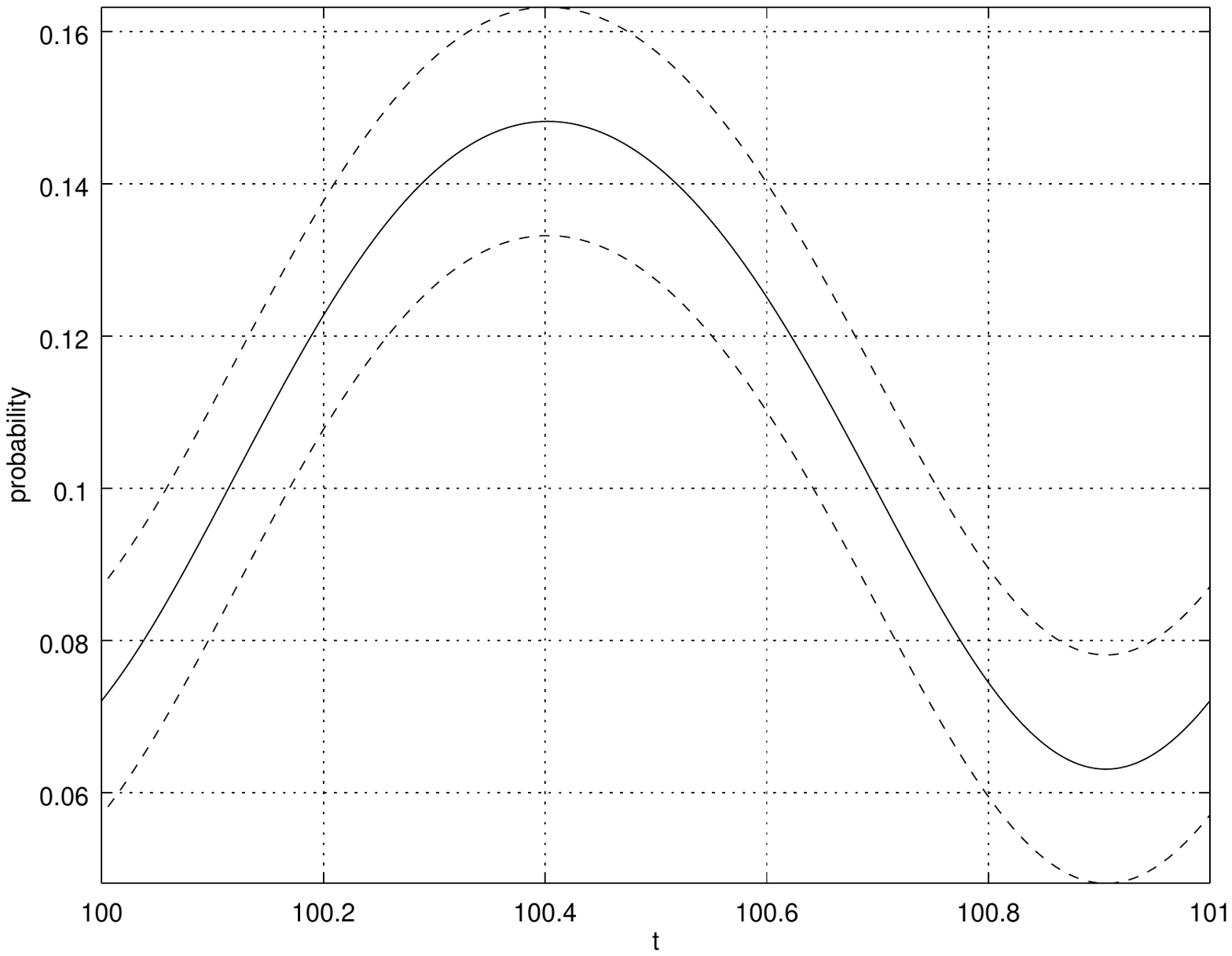}
\end{center}
\vspace{-5cm}\caption{Example 2. The perturbation bounds for the
'limit' probability $\Pr(X(t)=2)$, $t\in[100,101]$. }
\end{figure}

\end{document}